\documentclass[final]{siamltex}

\usepackage{amsmath,braket,amsfonts}
\usepackage[caption=false]{subfig}
\usepackage{graphicx,epstopdf}
\newtheorem{remark}{Remark}[section]
\newtheorem{example}{Example}[section]

\def\gamma{\mu}
\newcommand{\stl}
{\mathrel{\raise2pt\hbox{${\mathop<\limits_{\raise1pt\hbox{\mbox{$\sim$}}}}$}}}

\newcommand{\stlm}
{\mathrel{\raise2pt\hbox{${\mathop<\limits_{\raise1pt\hbox{\mbox{$\sim$}}}}_{\raise4.5pt\hbox{\begin{footnotesize}$\mu$\end{footnotesize}}}$}}}

\newcommand{\stg}
{\mathrel{\raise2pt\hbox{${\mathop>\limits_{\raise1pt\hbox{\mbox{$\sim$}}}}$}}}

\newcommand{\ste}
{\mathrel{\raise2pt\hbox{${\mathop=\limits_{\raise1pt\hbox{\mbox{$\sim$}}}}$}}}

\def\m{\mbox}
\def \no{\noindent}

\def \ee{\begin{equation}}
\def \e{\end{equation}}
\def\beq{\begin{eqnarray}}
\def\eq{\end{eqnarray}}
\def\beqx{\begin{eqnarray*}}
\def\eqx{\end{eqnarray*}}
\def\l{\label}  \def\no{\noindent}
\renewcommand{\theequation}{\arabic{section}.\arabic{equation}}

\def\omega{\alpha}

\def\ff{ {\mbox{\sc f}}}
\def\E{{\mbox{\sc e}}}
\def\vv{\mbox{\sc v}}
\def\o{\Omega}
\def\c{{\bf curl}}
\def\u{{\bf u}}
\def\v{{\bf v}}

\def\x{{\bf x}}
\def\ti{\times}
\def\la{{\bf \lambda}}

\def\p{\partial}
\def\q{\quad}

\def\w{{\bf w}}
\def\n{{\bf n}}
\def\f{{\bf f}}
\def\0{{\bf 0}}
\def\a{{\bf a}}
\def\b{{\bf b}}
\def\t{{\bf t}}
\def\r{{\bf r}}
\def\al{\alpha}
\def\be{\beta}
\def\T{{\mathcal T}_{h}}
\def\12{{1\over 2}}
\def\sl{\sup\limits}

\def\R{{\bf R}}

\def\chi{{\cal X}}

\def\cu{{\bf curl}}

\begin{document}

\title{Convergence of the Hiptmair-Xu Preconditioner for $H(\c)$-elliptic problems with Jump Coefficients (ii): Main Results}

\author{Qiya~Hu
\thanks{1. LSEC, ICMSEC, Academy of Mathematics and Systems Science, Chinese Academy of Sciences, Beijing
100190, China; 2. School of Mathematical Sciences, University of Chinese Academy of Sciences, Beijing 100049,
China (hqy@lsec.cc.ac.cn). This author was supported by the
Natural Science Foundation of China G12071469.}}


\maketitle
\begin{abstract}
This paper is the second of two articles, in which we aim to prove the convergence of the Hiptmair-Xu (HX) preconditioner (originally proposed by Hiptmair and Xu \cite{HX})
for $H(\c)$-elliptic boundary value problems with jump coefficients. In this paper, based on the auxiliary results obtained in our first article [Hu, SIAM J. Numer. Anal., 59(2021), pp.2500-2535],
we establish two new regular decompositions for low-dimensional edge finite element functions in three dimensions, which, under suitable assumptions on the distribution of the coefficients, are
stable up to a poly-logarithmic factor of the meshwidth with respect to weighted norms defined by the coefficients. Using these regular decompositions, we analyze the convergence of the
HX preconditioner for the case of strongly discontinuous coefficients. We show that the HX preconditioner is asymptotically optimal up to a poly-logarithmic
factor in the meshwidth and is not severely affected by large jumps of the coefficients across the interface between two neighboring subdomains.
\end{abstract}



\pagestyle{myheadings}
\thispagestyle{plain}
\markboth{}{}


\section{Introduction}\l{sec:introduction}
\setcounter{equation}{0}

On a bounded Lipschitz polyhedron $\Omega\subset \mathbb{R}^3$ consider the following ${\bf H}(\c,\Omega)$-elliptic boundary
value problem: \beq
\cu(\al\,\cu\,\u)+\be\u &=& \f\q \mbox{in}
\q\o,\cr \u\times\n&=&\0 \q \mbox{on} \q\o,
\label{1.1}
\eq
 where $\f\in (L^2(\Omega))^3$ is a given vector field,
$\al=\al(\x)$ and $\be=\be(\x)$ are spatially varying, scalar-valued, positive coefficients in $L^{\infty}(\Omega)$,
that are piecewise constant with respect to a partition of $\Omega$ into polyhedral Lipschitz sub-domains.

Generically, $\al$ and $\be$ will jump across interfaces between sub-domains. Variants of the boundary
value problem (\ref{1.1}) arise in various applications, for instance, as eddy current models in computational
electromagnetics \cite{Bos1998}. In this case, the coefficients characterize different homogeneous materials.

The most widely used Galerkin finite element discretization of (\ref{1.1}) is Ned$\acute{e}$l$\acute{e}$c's edge finite element
method \cite{s8}. Usually, the resulting sparse linear systems of equations are solved by some
preconditioned iterative method. As pointed out in \cite{DohrmannW2015} \cite{Tos} \cite{Wid05}, the construction of efficient preconditioners
for such linear systems is much more difficult than in the case of scalar second-order elliptic
boundary value problem discretized by means of Lagrangian finite elements.

Nevertheless, several efficient preconditioners for the edge finite element discretization of (\ref{1.1})
have been found, see, for example, \cite{Alo2} \cite{DohrmannW2015}
\cite{hipt} \cite{HX} \cite{HuShuZou2013} \cite{hz} \cite{HZ2} \cite{PasciakZhaos2002}
\cite{s10} \cite{Tos}. In particular, the HX preconditioner proposed
in \cite{HX} is very popular. The application of the HX preconditioner involves solving four discrete
scalar elliptic boundary value problems, which allows reusing existing codes.

Stable discrete Helmholtz decompositions or so-called regular decompositions are a key tool
in the convergence analysis of HX preconditioners. Using standard discrete regular decompositions,
their mesh-independent performance was shown for quasi-uniform sequences of meshes \cite{HX}. However, this theory cannot deal with large variations of the coefficients.
Although numerical results indicate that the HX preconditioner is not much affected by large jumps of the coefficients
\cite{Kol2009}, it remained an open problem whether the bounds from \cite{HX} for the condition numbers
achievable by HX preconditioning remain valid also in the case of coefficients $\al$ and $\be$ with large
jumps (refer to Subsection 7.3 of \cite{HX}).

This issue was first discussed in \cite{XuZ2011} for the case with two sub-domains, i.e., a simple transmission
problem. The main challenge was that stability estimates for standard discrete regular
decompositions had always relied on unweighted norms, which do not involve the coefficients
$\al$ and $\be$. The first important attempt towards a coefficient-aware discrete Helmholtz decomposition
was made in \cite{HuShuZou2013}, where a weighted discrete Helmholtz decomposition, which is almost
stable with respect to the weight function $\al=\be$, was constructed and studied. Moreover, using
this weighted Helmholtz decomposition in that paper the desired robust convergence result for
the preconditioner proposed in \cite{HZ2} was proved. Unfortunately, the weighted discrete Helmholtz
decomposition constructed in \cite{HuShuZou2013} cannot be applied to analyze the HX preconditioner for the
case with large jumps of the coefficients.

In the present paper, we introduce a new concept ``thorny vertex" associated with the distribution of jumps of the
coefficient $\al$ in order to reveal the essence of the problem of dealing with discontinuous coefficients.
Based on the auxiliary results obtained in \cite{Hu1-2017} and adapting some ideas
presented in \cite{HuShuZou2013}, we establish {\it new discrete regular decompositions}
for low-dimensional edge finite element functions in three dimensions, and show that the decompositions, on quasi-uniform
and shape-regular sequences of meshes and in weighted norms involving the coefficients $\al$ and
$\be$, are {\it nearly stable} with respect to the mesh size $h$ (i.e., stable up to a poly-logarithmic factor of $h$),
and {\it uniformly stable} with respect to large
jumps of the coefficients under suitable geometric constraints on the coefficients.

We would like to emphasize a key difference between the results obtained in this paper and \cite{HuShuZou2013}:
In this paper both $\al$ and $\be$ may have jumps and those may be independent of each other;
While \cite{HuShuZou2013} imposed restrictions on the behavior of the coefficients, those are no longer
needed for this work.

As an application of the proposed regular decompositions, we show that the PCG method with the HX preconditioner for solving
the considered system has a {\it nearly optimal} convergence rate. Nearly optimal means that the  convergence rate at worst
grows only poly-logarithmically like $O(|\log h|^q)$ with a positive integer $q$ as $h\rightarrow 0^+$,
and more importantly, does not suffer much from large jumps of the coefficients $\al$ and $\be$ across the interface between any two
neighboring subdomains.

Building on our first article \cite{Hu1-2017}, the present paper is organized as follows. In Section~2, we describe a domain decomposition based on
the distribution of coefficients and define some edge finite element subspaces. Section 3 introduces some new concepts and presents two new discrete regular
decompositions. Using these regular decompositions, we analyze the HX preconditioner for the case with
strongly discontinuous coefficients in Section 4. The new regular decompositions are constructed and analyzed in Section~5.

For the convenience of readers, we give a list of some important notations introduced and repeatedly used in this paper:
 \vskip 0.1in
 \begin{center}
\vskip 0.1in
\begin{tabular}{|c|c|} \hline
 symbol & position of its definition\\ \hline
 $\|\v\|_{H^{\ast}(\c\,,\o)}$&(\ref{norm1}) \\ \hline
 $\|\v\|_{H^1_{\ast}(\o)}$&(\ref{norm2}) \\ \hline
 ${\mathcal S}_k$ & (\ref{intersection}) \\ \hline
 $\Im^{\ast}_{\vv}$, $\Im^c_{\vv}$ & Definition \ref{definition4}  \\ \hline
 $n_{\ast}$ & (\ref{codimension}) \\ \hline
 $V^{\ast}_h(\Omega)$ & (\ref{co-space})    \\ \hline
 $\Im^{\ast}, \Im^c_{\ast}$ & (\ref{thorny-sets}) \\ \hline
 $\mathbb{P}(\mu)$ & Definition \ref{class} \\ \hline
 $\mathbb{P}_M(\mu)$ & Definition \ref{classM} \\ \hline
 $\stlm$ & in front of Theorem \ref{thm:main} \\ \hline
 \end{tabular}
 \end{center}

\section{Preliminaries}
This section introduces subdomain decompositions and some
fundamental finite element spaces.
\subsection{Sobolev spaces}
For an open and connected bounded domain ${\mathcal O}\subset \mathbb{R}^3$, let $H_0^1({\mathcal O})$ be the standard Sobolev space.
Define the Sobolev spaces
$$
H({\bf curl}; {\mathcal O}):=\{\v\in L^2({\mathcal O})^3; ~\cu\,\v\in L^2({\mathcal O})^3\}$$
and
$$ H_0({\bf curl}; {\mathcal O}):=\{\v\in
H({\bf curl}; {\mathcal O}); ~\v\times \n=0 ~\text{on} ~ \p{\mathcal O}\}.
$$

\subsection{Domain decomposition based on the distribution of
coefficients} \l{sec:discrete}

For precise descriptions of the new regular decompositions, we first decompose
the entire domain $\o$ into subdomains based on the surfaces of discontinuity of the coefficients
$\alpha(\x)$ and $\beta(\x)$ of (\ref{1.1}).

Associated with the coefficients $\al(\x)$ and $\be(\x)$, we assume that the entire domain $\o$ can be decomposed
into $N_0$ open Lipschitz polyhedral subdomains $\Omega_1, \Omega_2,\cdots,\Omega_{N_0}$ such that
$\bar\o=\cup_{k=1}^{N_0} \bar{\Omega}_k$ and the variations of the
coefficients $\al(\x)$ and $\be(\x)$ are not large in each
subdomain $\o_k$. Without loss of generality, we assume that, for
$k=1,2, \ldots,N_0$,
\begin{equation}
\label{eq:coeffb}
\omega(\x)=\omega_k~~\mbox{and}~~\beta(\x)=\beta_k,
\q\forall\,\x\in \Omega_k,
\end{equation}
where each $\al_k$ or $\be_k$ is a positive constant. Such a decomposition is possible in many applications
when $\o$ contains several media. The results established in this paper still hold if the coefficients $\al(x)$ and $\be(x)$ are
not constants but have small variations in some subdomains. Note that a subdomain $\Omega_k$ may be a non-convex polyhedron, which is a union of several convex polyhedra.
In this sense, our assumption is not restrictive and does cover many practical cases.

The subdomains $\{\Omega_k\}_{k=1}^{N_0}$ are of different nature from those in the context of the standard domain decomposition methods: $\{\Omega_k\}_{k=1}^{N_0}$
is decomposed based only on the distribution of jumps of the coefficients $\alpha(x)$ and $\beta(x)$ (so $N_0$ is a fixed integer, and the size of each
$\Omega_k$ is $O(1)$).

\subsection{Edge and nodal element spaces}

Let $\Omega$ be equipped with a family of tetrahedral triangulations ${\mathcal T}_h$ with the mesh size $h$.
We assume that the families of triangulations ${\mathcal T}_h$ with $h\rightarrow 0^+$ are shape regular and quasi-uniform (refer to \cite{Hu1-2017}).
The Nedelec edge element space, of the lowest order, is a subspace of piecewise linear
polynomials defined on ${\mathcal T}_h$ (refer to \cite{brezzi} and \cite{s8}): 
$$ \l{eq:nedelec}
V_h(\o):=\Big\{\v\in {H}_0({\bf curl};\o); ~\v|_K\in R(K), ~\forall K\in\T
 \Big\},
$$
where $R(K)$ is a subset of all linear polynomials on the element
$K$ of the form:
$$ R(K)=\Big\{\a+\b\ti\x;~ \a,\b\in {\bf R}^3, ~\x\in K\Big\} \,. $$

Let ${\mathcal E}_h$ denote the set of edges of ${\mathcal T}_h$. Then, for any $\v\in V_h(\o)$, its tangential components
are continuous on all edges $e\in{\mathcal E}_h$, and $\v$ is uniquely determined by its moments on all fine
edges $e$:
$$ M_h(\v)=\Big\{\la_e(\v)=\int_e\v\cdot\t_e ds; ~ e\in {\mathcal E}_h\Big\} $$
where $\t_e$ denotes the unit vector on edge $e$. We write $e$ to denote any edge or union of edges, either from
an element $K\in\T$ or from a subdomain. For a vector-valued
function $\v$ with appropriate smoothness, we introduce its edge
element ``projection" ${\bf r}_h\v$ (refer to \cite{Hu1-2017}) such that ${\bf r}_h\v\in
V_h(\Omega)$, and ${\bf r}_h\v$ and $\v$ have the same moments as given in $M_h(\v)$.

As we will see, the edge element analysis also frequently involves the
nodal element space. For this purpose, we introduce $Z_h(\o)$ to be
the standard continuous piecewise linear finite element space associated with the
triangulation ${\mathcal T}_h$.

%

Let $G$ be a subdomain of $\o$, where $G$ is a union of elements in ${\mathcal T}_h$. Associated with $G$, we respectively use $V_h(G)$ and $Z_h(G)$ to
denote the natural restriction of $V_h(\o)$ and $Z_h(\o)$ on $G$.

\subsection{Weighted norms}

Let $\al$ and $\be$ denote the coefficients appearing in (\ref{1.1}). For $H(\c)$ functions, we define
$$
\|\v\|_{L^2_{\omega}(\o)}:=(\sum\limits_{r=1}^{N_0}\omega_r\|\v\|^2_{0,\o_r})^{\12},\quad\v\in H(\c;~\Omega)
$$
and
\ee
\|\v\|_{H^{\ast}(\c\,,\o)}:=(\sum\limits_{r=1}^{N_0}\omega_r\|\c~\v\|^2_{0,\o_r}+\be_r\|\v\|^2_{0,\o_r})^{\12},\quad\v\in H(\c;~\Omega).\label{norm1}
\e
For $H^1$ functions, we define
$$
\|p\|_{H^1_{\beta}(\o)}:=(\sum\limits_{r=1}^{N_0}\beta_r\|\nabla p\|^2_{0,\o_r}+\beta_r\|p\|^2_{0,\o_r})^{\12},\quad p\in H^1(\Omega)
$$
and
\ee
\|\v\|_{H^1_{\ast}(\o)}:=(\sum\limits_{r=1}^{N_0}(\omega_r|\v|^2_{1,\o_r}+\be_r\|\v\|^2_{0,\o_r}))^{\12},\quad\v\in (H^1(\Omega))^3.\label{norm2}
\e
The stability estimates derived in this paper will be based on the above weighted norms.


\section{Main results}\label{sec:helmholtz}
\setcounter{equation}{0}
In this section, we present new discrete regular decompositions associated with the distribution of discontinuous coefficients
such that the decompositions are nearly stable with respect to the weighted norms defined in Subsection 2.4.

It is often a knotty problem to construct and analyze preconditioners for partial differential equations with jump coefficients.
If there are {\it internal cross-points} in the distribution of jump coefficients, multilevel preconditioners
for second-order elliptic problems are not uniformly convergent with respect to the jumps of the coefficients,
unless the distribution of coefficients satisfies some additional assumptions, for example, the {\it quasi-monotonicity assumption} (see \cite{DSM1996} and \cite{Pet2002}),
which is a relatively weak assumption.
In this section, we try to establish discrete regular decompositions for general cases without the {\it quasi-monotonicity assumption}.
To this end, we need to introduce some new concepts and notations.

\subsection{Upper intersection sets and thorny vertices}

In this subsection, we introduce two new concepts {\it upper intersection set} and {\it thorny vertex}.

In the following, we always use $\vv$ to denote a vertex of one polyhedron obtained from the domain decomposition described in Subsection 2.2.
Let ${\mathcal N}(\Omega)$ (resp. ${\mathcal N}(\partial\Omega)$) denote the set of vertices $\vv$ in $\Omega$ (resp. on $\partial\Omega$).
For a vertex $\vv\in{\mathcal N}(\Omega)\cup{\mathcal N}(\partial\Omega)$, let $\Im_{\vv}$ denote the set of all
polyhedra $\Omega_k$ that contain $\vv$ as one of their vertices.

Let $\Omega_1,\cdots,\Omega_{N_0}$ be the subdomains defined in Subsection 2.2. For a subdomain $\o_k$ ($1\le k\le N_0$),  let ${\mathcal S}_k$ be the union of the
(closed) intersection sets of $\partial\Omega_k$ with $\partial\o$ or $\partial\Omega_l$ ($l\not=k$) that satisfy $\alpha_l\geq\al_k$, namely,
\ee
{\mathcal S}_k:=\bigcup_{\alpha_l\geq\alpha_k,l\not=k}(\partial\Omega_k\cap\partial\Omega_l)\bigcup(\partial\Omega_k\cap \partial\o).\label{intersection}
\e
When the maximum $\mbox{max}_{1\leq r\leq N_0}\al_r$ is achieved on $\Omega_k$, we have ${\mathcal S}_k=\partial\Omega_k\cap \partial\o$, which may be the empty set.
In general the set ${\mathcal S}_k$ is a union of some faces, edges and vertices, which may be non-connected. For convenience, if a vertex in ${\mathcal S}_k$
does not intersect with any face or edge in ${\mathcal S}_k$, then it is called an ``isolated vertex" in ${\mathcal S}_k$.
Similarly, we can define an ``isolated edge" in ${\mathcal S}_k$.
\begin{definition}\label{definition-new} The set ${\mathcal S}_k$ defined by (\ref{intersection}) is called the {\it upper intersection set} of the subdomain $\Omega_k$.
\end{definition}

The structure of the {\it upper intersection set} ${\mathcal S}_k$ can describe the distribution of the coefficients corresponding to all the subdomains adjacent to $\Omega_k$.

\begin{definition}\label{definition3} A vertex $\vv\in{\mathcal N}(\Omega)\cup{\mathcal N}(\partial\Omega)$ is called a {\it thorny vertex} with respect to the distribution of the
coefficients $\al_1,\al_2,\cdots,\al_{N_0}$ if there exists a polyhedron $\Omega_k\in\Im_{\vv}$ such that $\vv$ is an isolated vertex in its {\it upper intersection set} ${\mathcal S}_k$.
\end{definition}

A {\it thorny vertex} $\vv$ can also be defined as follows. There exists a subdomain $\Omega_k\in\Im_{\vv}$ such that: (i) if $\vv\in{\mathcal N}(\Omega)$, for any other polyhedron $\Omega_{r'}$
that belongs to $\Im_{\vv}$ and corresponds to the coefficient satisfying $\al_{r'}\geq\al_r$, the intersection of $\bar{\Omega}_r$ with $\bar{\Omega}_{r'}$ is just the vertex $\vv$, i.e., $\bar{\Omega}_r\cap\bar{\Omega}_{r'}=\vv$; (ii) if $\vv\in{\mathcal N}(\partial\Omega)$,
the intersection of $\bar{\Omega}_r$ with $\partial\Omega$ is just the vertex $\vv$, and the local maximum of $\al(\x)$ on the union of all the polyhedra in $\Im_{\vv}$ is achieved on $\bar{\Omega}_r$.

A {\it thorny vertex} $\vv$ in ${\mathcal N}(\Omega)$ (resp. ${\mathcal N}(\partial\Omega)$) is called
an {\it internal thorny vertex} (resp. a {\it boundary thorny vertex}).

A typical example of {\it internal thorny vertex} is the central vertex in the ``checkerboard" domain described below.

\begin{example}\label{example3} Set $\Omega=[0,1]^3$, and let $\Omega$ be decomposed into a union of four polyhedra $\{\Omega_k\}_{k=1}^4$, where (see Figure 4)
$$\Omega_1=[0,{1\over 2}]^3,~~\Omega_2=[{1\over 2},1]^3,~~\Omega_3=([0,1]^2\times[0,{1\over 2}])\backslash \Omega_1,~~\Omega_4=([0,1]^2\times[{1\over 2},1])\backslash \Omega_1.$$
Define $\al(\x)=1$ on $\Omega_1\cup\Omega_2$, and $\al(\x)=\varepsilon<1$ on $\Omega_3\cup\Omega_4$, i.e, $\al_1=\al_2=1$ and $\al_3=\al_4<1$.
\end{example}
\vskip 2mm
\begin{center}
\includegraphics[width=8cm,height=6cm]{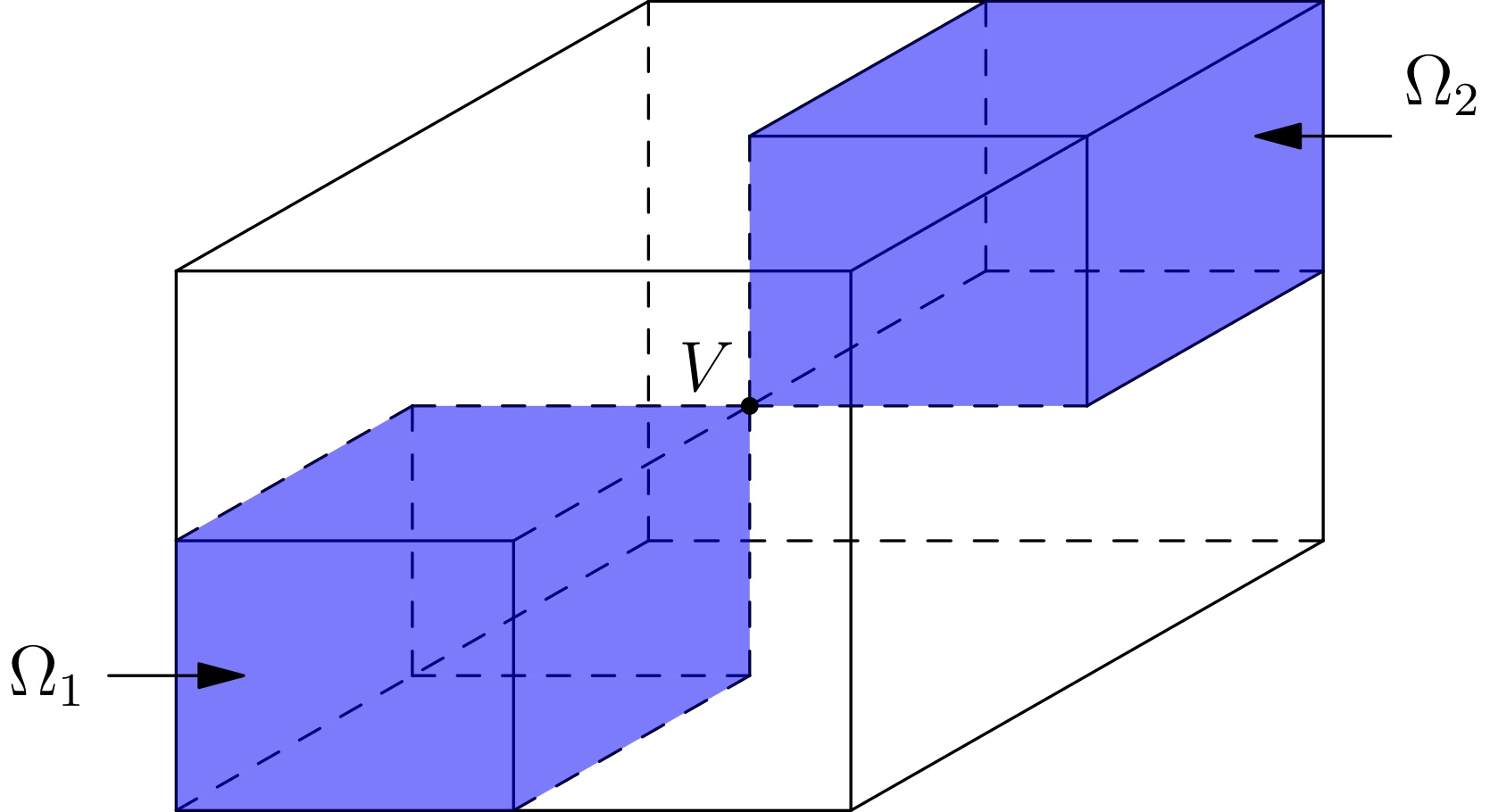}
\vskip 0.2in
\centerline{Figure 1: A ``checkerboard" domain $\Omega$: the vertex $\vv$ is an {\it internal thorny vertex}.}
\end{center}
\vskip 2mm

The domain shown in Figure 1 is called {\bf ``checkerboard" domain}. In this example, the vertex $\vv$ at the center of $\Omega$ is an {\it internal thorny vertex} because $\vv$ is an isolated
vertex in the {\it upper intersection set} ${\mathcal S}_2$.

Let us give another example to explain Definition \ref{definition-new} and Definition \ref{definition3}.

\begin{center}
\includegraphics[width=6.5cm,height=4.3cm]{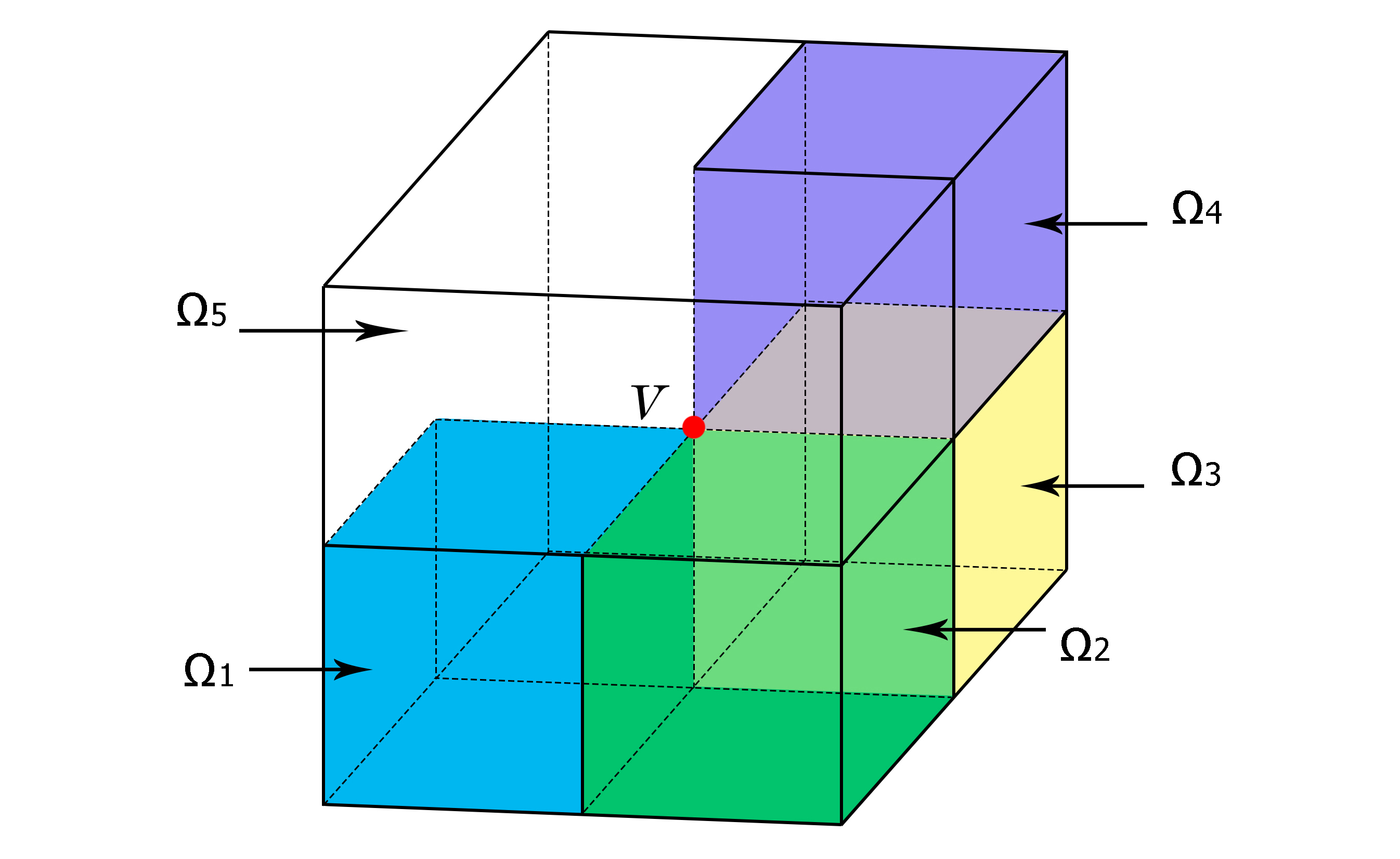}\quad

\centerline{Figure 2: {\it upper intersection set} and {\it thorny vertex}.
}
\end{center}
\vskip 2mm
\begin{example}\label{example1} Consider the subdomains shown in Figure 2.

Case (i): $\al_1>\al_3>\al_2>\al_4>\al_5$. The {\it upper intersection set} ${\mathcal S}_3$ contains an isolated edge $\partial\Omega_3\cap\partial\Omega_1$ but the vertex $\vv$ is not a
{\it thorny vertex}.

Case (ii): $\al_1>\al_4>\al_2>\al_3>\al_5$. The vertex $\vv$ is a {\it internal thorny vertex} because it is an isolated
vertex in the {\it upper intersection set} ${\mathcal S}_4$.

\end{example}

A {\it thorny vertex} $\vv$ can be roughly understood as: there exist two polyhedra such that their intersection is just $\vv$ and they correspond to larger coefficients than
any other polyhedron containing $\vv$ as its vertex.


\begin{definition}\label{definition4} Let $\vv$ be a {\it thorny vertex}. Two subsets $\Im^{\ast}_{\vv}$ and $\Im^c_{\vv}$ of $\Im_{\vv}$ are described as follows:
a polyhedron $\Omega_r\in \Im^{\ast}_{\vv}$, if and only if each polyhedron $\Omega_{r'}$ belonging to $\Im_{\vv}$ and satisfying $\al_{r'}\geq\al_r$ has the property $\bar{\Omega}_r\cap\bar{\Omega}_{r'}=\vv$;
a polyhedron $\Omega_r\in \Im^c_{\vv}$ if and only if there is at least one polyhedron $\Omega_{r'}\in\Im_{\vv}$ such that $\al_{r'}\geq\al_r$ and the intersection of $\Omega_r$ and $\Omega_{r'}$
is a face or an edge containing $\vv$.
\end{definition}

According to Definition \ref{definition4}, for Case (ii) of Example \ref{example1} we have $\Im^{\ast}_{\vv}=\{\Omega_1,\Omega_4\}$ and $\Im^c_{\vv}=\{\Omega_2,\Omega_3,\Omega_5\}$ (see Figure 2).

For any two polyhedrons $\Omega_r,~\Omega_l\in \Im^{\ast}_{\vv}$, we have $\bar{\Omega}_r\cap\bar{\Omega}_{l}=\vv$.

\begin{remark} The set $\Im_{\vv}$ can be decomposed into a union of two disjoint sets $\Im^{\ast}_{\vv}$ and $\Im^c_{\vv}$. Moreover, neither $\Im^{\ast}_{\vv}$ nor $\Im^c_{\vv}$
is the empty set. In fact, $\Im^{\ast}_{\vv}$ contains at least $\Omega_k$ mentioned in Definition \ref{definition3} and the polyhedron achieving the locally maximal coefficient.
Note that $\Im^{\ast}_{\vv}$ is a real set of $\Im_{\vv}$ (otherwise, $\vv$ is not a {\it thorny vertex}), so $\Im^c_{\vv}$ is not the empty set.
\end{remark}

The following result follows directly from the definitions of $\Im^{\ast}_{\vv}$ and $\Im^c_{\vv}$.\\
{\bf Proposition 3.1}. Let $\vv$ be a {\it thorny vertex} and $\Omega_r\in\Im^{\ast}_{\vv}$. Then, for each polyhedron $\Omega_l$ that belongs to $\Im_{\vv}$ and has a common
face or a common edge with $\Omega_r$, we have $\Omega_l\in\Im^c_{\vv}$ and $\al_l<\al_r$.

As we will see, the introduction of {\it thorny vertices} helps us to build a weighted-stable regular decomposition of edge element functions for more complicated
coefficient distributions.

\subsection{A discrete regular decomposition}

In this subsection, we describe a regular decomposition for the case with {\it thorny vertices}.

Let ${\mathcal V}_{\ast}$ denote the set of all {\it thorny vertices} generated by the domain decomposition and the distribution of coefficients.
For $\vv\in{\mathcal V}_{\ast}$, we use $n_{\vv}$ to denote the number of subdomains contained in $\Im^{\ast}_{\vv}$.
Let ${\mathcal V}^{in}_{\ast}$ (resp. ${\mathcal V}^{\partial}_{\ast}$) denote the set of {\it internal thorny vertices} (resp. {\it boundary thorny vertices}). Define
\ee
 n_{\ast}:=\sum\limits_{\vv\in{\mathcal V}^{in}_{\ast}}(n_{\vv}-1)+\sum\limits_{\vv\in{\mathcal V}^{\partial}_{\ast}}n_{\vv}. \label{codimension}
\e
For convenience, the number $n_{\ast}$ is called the {\it multiplicity of thorny vertices}. This value reflects the number of {\it thorny vertices} and the
distribution of coefficients on the polyhedron subdomains containing a {\it thorny vertex}.

As we will see in {\bf Proposition 3.3}, a function $\v_h\in V_h(\Omega)$ does not automatically admit a stable discrete regular decomposition when ${\mathcal V}_{\ast}\not=\emptyset$,
unless $\v_h$ meets an additional condition associated with each {\it thorny vertex}.
Such a condition can be formally written as ${\mathcal F}_{\vv}\v_h=0$, where ${\mathcal F}_{\vv}$ denotes a functional associated with the {\it thorny vertex} $\vv$.
The functional ${\mathcal F}_{\vv}$ can be defined in different manners and its different definitions do not affect the application of
the main results (see the next section). In this subsection, we present a special definition of ${\mathcal F}_{\vv}$, as in Theorems 3.3 and 3.4
in Section 3 of \cite{Hu1-2017}. To define such a functional, we first introduce some notation.

For $\vv_i\in{\mathcal V}_{\ast}$, let $\{\Omega_{i_r}\}_{r=1}^{m_i}$ denote all the subdomains in $\Im^{\ast}_{\vv_i}$. For $\Omega_{i_r}$ ($r=1,\cdots,m_i$), we choose a face $\ff_{i_r}$
of $\Omega_{i_r}$ such that $\vv_i$ is a vertex of $\ff_{i_r}$, and further choose an edge $\E_{i_r}$ of $\ff_{i_r}$ (see Figure 3 for the simplest case with $m_i=2$).
When $\partial\ff_{i_r}\cap\partial\Omega\not=\emptyset$ (so it is a union of edges of $\Omega_{i_r}$), we choose $\E_{i_r}\subset\partial\ff_{i_r}\cap\partial\Omega$.

\begin{center}
\includegraphics[width=9cm,height=6cm]{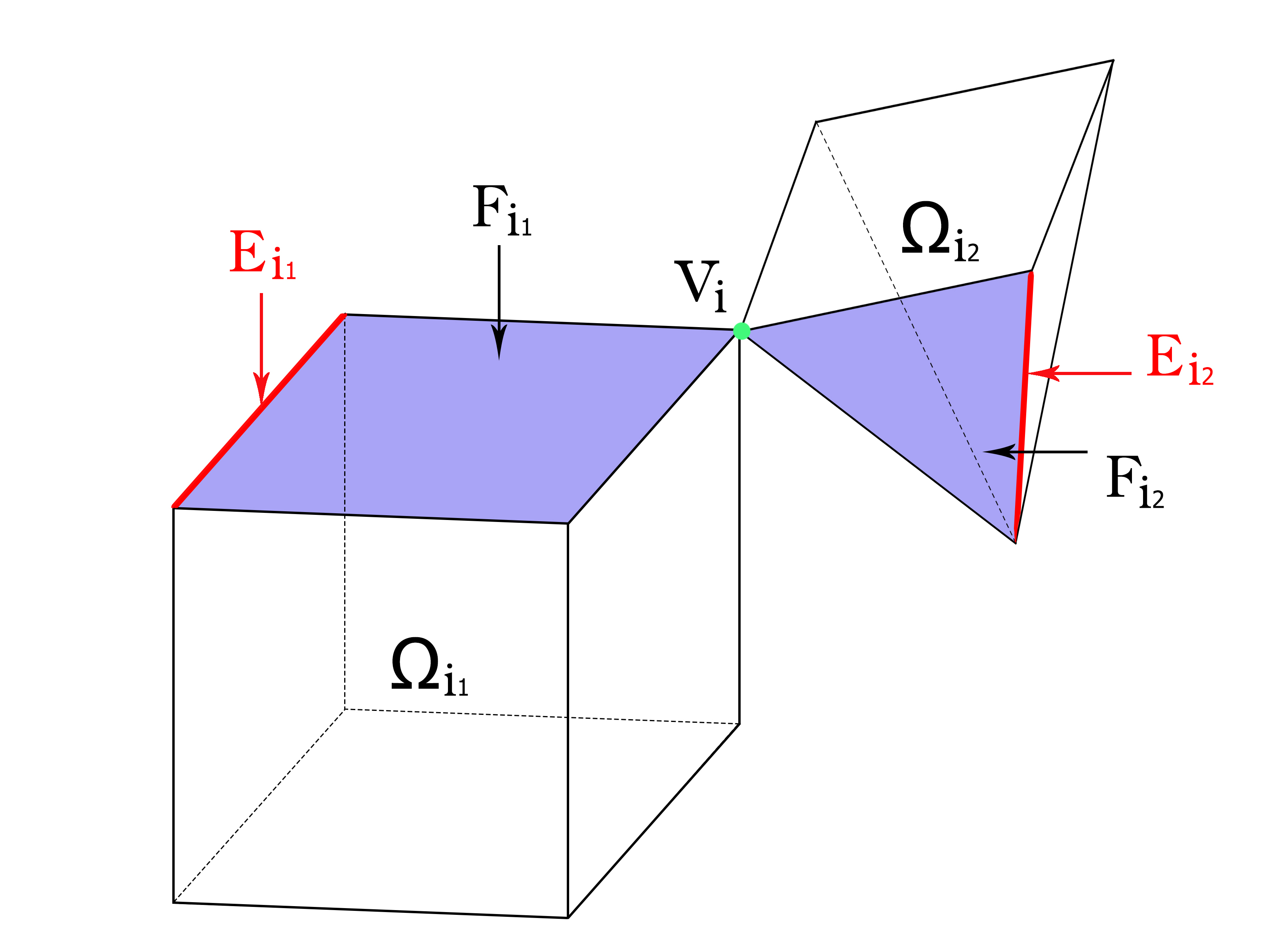}
\vskip 0.1in
\centerline{Figure 3: Two polyhedrons that share an {\it internal thorny vertex}}
\end{center}
\vskip 2mm

Let $s$ be the arc-length variable along $\partial\ff_{i_r}$, taking values from $0$ to $l_{i_r}$, where $l_{i_r}$ is the total length of $\partial\ff_{i_r}$ and $t=0$ (and $t=l_{i_r}$)
corresponds to the vertex $\vv_i$. As in \cite{Tos}, for an edge finite element function $\v_h$, we define ($r=1,\cdots,m_i$)
$$
C^{\partial\ff_{i_r}}_{\v_h}=\frac{1}{l_{i_r}}\int_0^{l_{i_r}}(\v_h\cdot\t_{\partial\ff_{i_r}})(s)\,ds $$
and
$$
\phi^{\partial\ff_{i_r}}_{\v_h}(t)
=\int_0^{t} (\v_h\cdot\t_{\partial\ff_{i_r}}-C^{\partial\ff_{i_r}}_{\v_h})(s)ds+c_{i_r}\,,
~\forall\,t\in [0, l_{i_r}]\,,
$$
where the constant $c_{i_r}$ is chosen such that $\phi^{\partial\ff_{i_r}}_{\v_h}$ has an average value of zero on $\E_{i_r}$. Then,
the tangential component $\v_h\cdot\t_{\partial\ff_{i_r}}$ of $\v_h$ admits the decomposition (see \cite{Tos})
$$ (\v_h\cdot\t_{\partial\ff_{i_r}})(t)=C^{\partial\ff_{i_r}}_{\v_h}+(\phi^{\partial\ff_{i_r}}_{\v_h})'(t),\quad t\in [0, l_{i_r}]\quad (r=1,\cdots,m_i). $$

Define the functional
$$ {\mathcal F}_{i,r}\v_h:=\phi^{\partial\ff_{i_1}}_{\v_h}(l_{i_1})-\phi^{\partial\ff_{i_r}}_{v_h}(l_{i_r}),\quad\quad r=2,\cdots,m_i \quad\quad (\vv_i\in{\mathcal V}^{in}_{\ast}) $$
or
$$ {\mathcal F}_{i,r}\v_h:=\phi^{\partial\ff_{i_r}}_{v_h}(l_{i_r}),\quad\quad r=1,\cdots,m_i\quad\quad (\vv_i\in{\mathcal V}^{\partial}_{\ast}). $$
With the help of the functionals $\{{\mathcal F}_{i,r}\}$, we can define a subspace
\ee
V^{\ast}_h(\Omega):=\{\v_h\in V_h(\Omega):~{\mathcal F}_{i,r}\v_h=0~~\mbox{for~each}~r~\mbox{and}~~i=1,\cdots,n_{\ast}\}.\label{co-space}
\e

Note that the condition ${\mathcal F}_{i,r}\v_h=0$ implies that
$$ \phi^{\partial\ff_{i_r}}_{v_h}(l_{i_r})=\phi^{\partial\ff_{i_1}}_{\v_h}(l_{i_1}),~~~r=2,\cdots,m_i\quad (\vv_i\in{\mathcal V}^{in}_{\ast}) $$
and $t=l_{i_r}$ corresponds to the {\it thorny vertex} $\vv_i$ for $\Omega_{i_r}$, the condition ${\mathcal F}_{i,r}\v_h=0$ means that $\v_h$ should have some ``continuity" at $\vv_i\in{\mathcal V}^{in}_{\ast}$.
Of course, when there is no {\it thorny vertex}, we have $V^{\ast}_h(\Omega)=V_h(\Omega)$.
It is clear that $\mbox{dim}(V^{\ast}_h(\Omega))=\mbox{dim}(V_h(\Omega))-n_{\ast}$,
so the {\it multiplicity of thorny vertices} is just the {\it codimension} of the space $V^{\ast}_h(\Omega)$.

For convenience, we introduce the following two sets:
\ee
\Im^{\ast}=\cup_{\vv\in {\mathcal V}_{\ast}}\Im^{\ast}_{\vv},\quad\quad\Im^c_{\ast}=\{\Omega_k\}_{k=1}^{N_0}\backslash\Im^{\ast}.\label{thorny-sets}
\e
{\bf Proposition 3.2}. For each $\o_r\in\Im^c_{\ast}$, the {\it upper intersection set} ${\mathcal S}_r$ does not contain any isolated vertex, i.e., ${\mathcal S}_r$ is a union of faces and edges only.\\
\no{\it Proof}. Let $\o_r\in\Im^c_{\ast}$. Assume that $\vv$ is an isolated vertex in ${\mathcal S}_r$. By the definition of ${\mathcal S}_r$, there is a polyhedron $\Omega_l$ satisfying $\alpha_l\geq\al_r$
such that $\bar{\Omega}_l\cap\bar{\Omega}_r=\vv$, which means that $\vv$ is a {\it thorny vertex} by Definition \ref{definition3}. Moreover, for any polyhedron $\Omega_{r'}$ belonging to $\Im_{\vv}$ and satisfying
$\al_{r'}\geq \al_r$, we have $\bar{\Omega}_{r'}\cap\bar{\Omega}_r=\vv$ (otherwise, $\vv$ is not an isolated vertex). Thus $\Omega_r\in\Im^{\ast}_{\vv}\subset\Im^{\ast}$ by Definition \ref{definition4},
which contradicts the condition that $\Omega_r\in\Im^c_{\ast}$. Therefore, ${\mathcal S}_r$ is a union of faces and edges only. \hfill $\Box$

As in \cite{Hu1-2017}, a connected union $\Gamma$ of some faces is called {\it a connected Lipschitz union of some faces} if there is no isolated vertex in $\Gamma$
(see Figure 4), a more detailed explanation was given in Subsection 2.3 of \cite{Hu1-2017}.
A union $\hat{\Gamma}$ of some faces is called a {\it union of connected Lipschitz unions of faces} if
$\hat{\Gamma}$ is a (possibly non-connected) union of $\Gamma_1,\cdots,\Gamma_J$, with $\Gamma_j$ being a {\it connected Lipschitz union of some faces} ($j=1,\cdots,J$).

\begin{center}
\includegraphics[width=5cm,height=5cm]{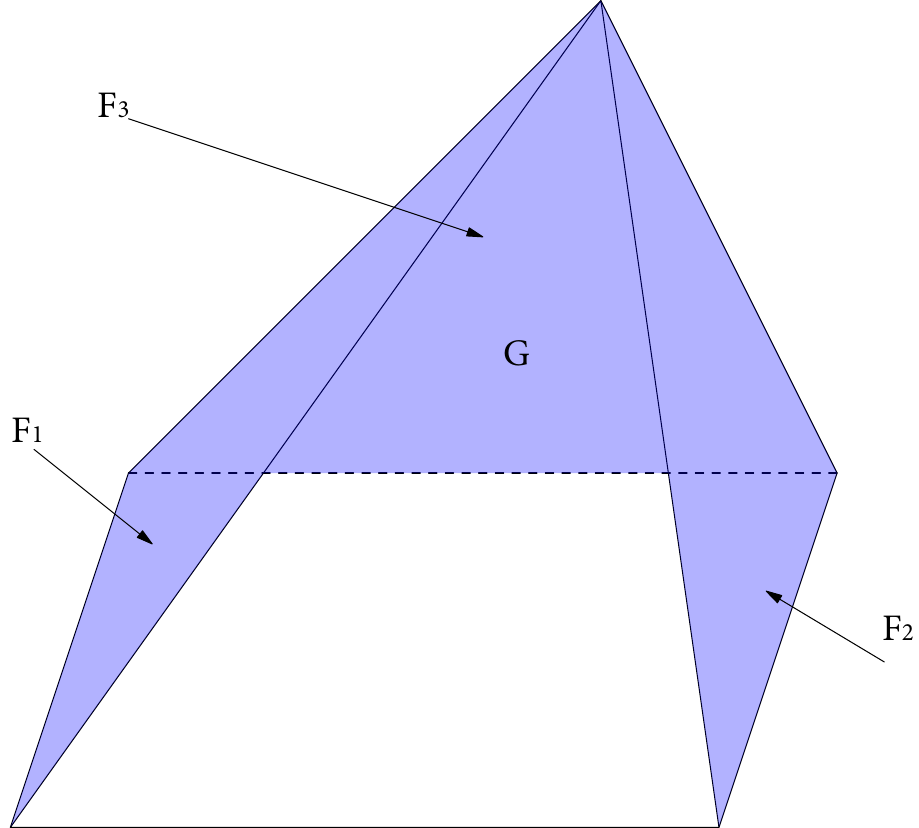}\quad
\includegraphics[width=5cm,height=5cm]{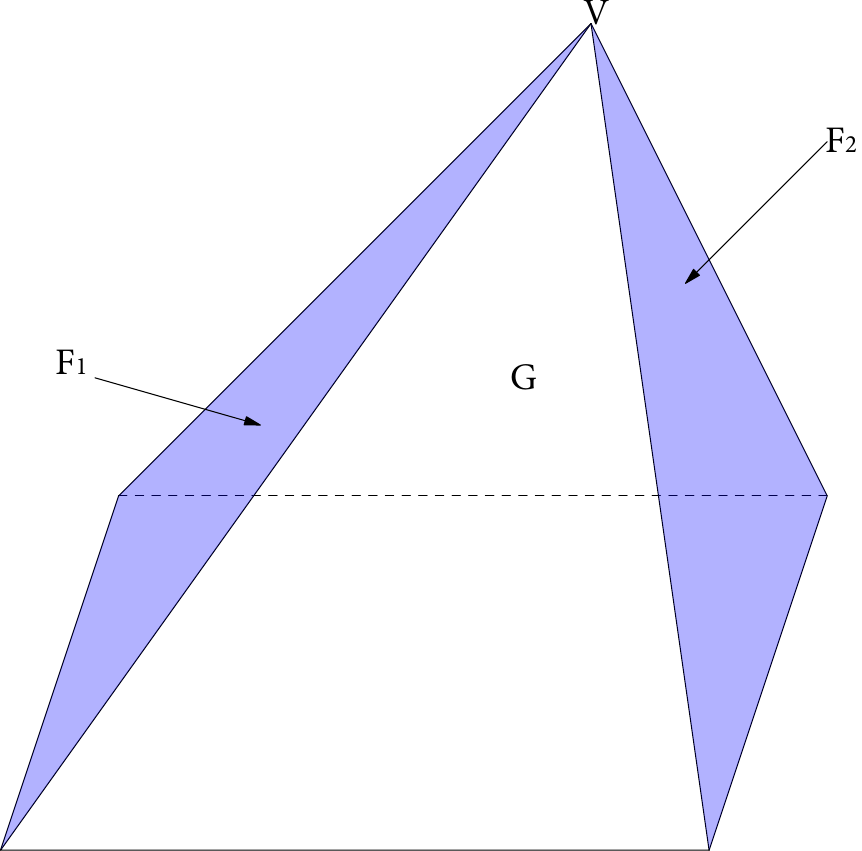}
\end{center}
\vskip 0.1in
\begin{minipage}{0.90\linewidth}
Figure 4. $\ff_1$ and $\ff_2$ are two opposite lateral faces of the four-sided pyramid $G$, and $\ff_3$
denotes another lateral face that has a common edge with each of $\ff_1$ and $\ff_2$. The
set $\Gamma=\ff_1\cup\ff_2\cup\ff_3$ is a {\it connected Lipschitz union of three faces}~(left), but the
set $\Gamma=\ff_1\cup\ff_2$ is a {\it connected non-Lipschitz union of two faces}~(right), where the
point $\vv$ is an isolated vertex of $\Gamma$.
\end{minipage}
\vskip 0.1in
We define a function $h\in \mathbb{R}^+\longmapsto \rho(h)$ as follows: $\rho(h)=1$ if, for all the subdomains $\Omega_k$ in $\Im^c_{\ast}$, the {\it upper intersection set} ${\mathcal S}_k$
is a union of connected Lipschitz unions of some faces of $\Omega_k$; otherwise, $\rho(h)=\log(1/h)$.

In applications, the coefficients $\{\al_k\}$ and $\{\beta_k\}$ may be dependent on or independent of each other (for the latter case, the ratio $\al_k/\beta_k$ may be large or small).
Then all the coefficients $\{\al_k\}$ and $\{\beta_k\}$ have three possible relations: (1)  the ratio $\be_k/\al_k$ is nether large nor small;
(2) the ratio $\be_k/\al_k$ is not small but may be large; (3) the ratio $\be_k/\al_k$ is not large but may be small.
As we will see, for different situations, we need to use different stability estimates in the regular decompositions developed in \cite{Hu1-2017}. For example, for the third case,
we have to use the stability of a vector-valued $H^1$ function with respect to the $\c$ semi-norm (a complete norm-controlled stability of the vector-valued $H^1$ function is not practical for this case).
To cover all three situations, we have to carefully investigate the coefficients.

\begin{definition}\label{class} Let $\mu\geq 1$ be a given constant (independent of
the mesh size $h$). We say that the coefficients $\{\al_k\}$ and $\{\be_k\}$ belong to the problem class $\mathbb{P}(\gamma)$ if
the coefficients $\{\al_k\}$ and $\{\be_k\}$ satisfy three conditions:

{\bf Condition A}. The coefficients on any two neighboring subdomains $\Omega_i$ and $\Omega_j$ have the property
\ee
\be_i\leq \gamma \be_j~~~\mbox{when}~~\al_i\leq\al_j. \label{new2.0001}
\e

{\bf Condition B}. For every $\Omega_k\in \Im^c_{\ast}$, on which the coefficients do not satisfy the relation $\gamma^{-1}\al_k\leq \be_k\leq\gamma\al_k$,
the {\it upper intersection set} ${\mathcal S}_k$ (if it is non-empty) meets one of the following two requirements:

(i) if $\be_k>\mu\al_k$, then ${\mathcal S}_k$ is a (possibly non-connected) union of some {\it faces} of $\Omega_k$;

(ii) if $\be_k<\mu^{-1}\al_k$, then ${\mathcal S}_k$ is a {\it connected} union of some faces and edges of $\Omega_k$.

{\bf Condition C}. The coefficients $\al_k$ and $\be_k$ for $\Omega_k\in \Im^{\ast}$ have the relation:
\begin{equation}
\label{new2.1}
\beta_k\leq\mu\alpha_k,~~~\forall\Omega_k\in \Im^{\ast}.
\end{equation}

\end{definition}

In the above {\bf Condition B}, when the coefficients on $\Omega_k$ satisfy the relation $\gamma^{-1}\al_k\leq \be_k\leq\gamma\al_k$, the set ${\mathcal S}_k$ need not meet a stricter requirement than that
${\mathcal S}_k$ is a (possibly non-connected) union of some faces and edges of $\Omega_k$, which is naturally met by {\bf Proposition 3.2}.

\begin{remark} It can be seen from Definition \ref{class} that $\mathbb{P}(\mu_1)\subset \mathbb{P}(\mu_2)$ when $\mu_1\leq \mu_2$. The smallest class $\mathbb{P}(1)$ ($\mu=1$) contains
all the coefficients $\{\al_k\}$ and $\{\be_k\}$ with $\be=\al$. Assume that $c_0\al\leq\be\leq C_0\al$ for two given positive constants $c_0$ and $C_0$, which is
a common assumption in existing works \cite{HX,HuShuZou2013}.
Then the coefficients $\{\al_k\}$ and $\{\beta_k\}$ belong to the class $\mathbb{P}(\mu)$ with $\mu=\max\{C_0,1/c_0\}$ (without additional assumption on ${\mathcal S}_k$).
\end{remark}

From now on, we shall always use $C$ to denote a generic constant only depending on the sub-domains
$\Omega_k$, shape-regularity and quasi-uniformity of the mesh, but independent of the meshwidth $h$ and the coefficients $\al$ and $\be$.
For convenience, we introduce two new notations $\stl$ and $\stlm$. For any two non-negative quantities $x$ and $y,$  $x\stl y$
means that $x\leq C y$ for some constant $C$ described above, and $x\stlm y$
means that $x\leq C(\mu) y$ for some constant $C(\mu)$ only depending on the sub-domains
$\Omega_k$, shape-regularity and quasi-uniformity of the mesh and the parameter $\mu$ in the problem classes $\mathbb{P}(\mu)$ and
$\mathbb{P}_M(\mu)$ defined later, but independent of the meshwidth $h$ and
the specific values and the detailed structure of the coefficients $\al$ and $\be$.

The following theorem presents a general result for a weighted discrete regular decomposition.

\begin{theorem}\l{thm:main} Assume that the coefficients $\{\al_k\}$ and $\{\be_k\}$ belong to the problem class $\mathbb{P}(\gamma)$.
Then, any function $\v_h\in V^{\ast}_h(\o)$ admits a decomposition of the form
\ee
  \v_h=\nabla p_h+\r_h\w_h+{\bf R}_h \label{newnew3.01}\e for some
  $p_h\in Z_h(\o)$ and $\w_h\in (Z_h(\o))^3$ and ${\bf R}_h\in
  V_h(\o)$. Moreover, we have the estimates
  \ee \|\beta^{\12}\nabla p_h\|_{0,\o}\stlm
\log^2(1/h)\rho^{m_0}(h)\|\v_h\|_{H^{\ast}(\c,~\o)}\label{stab3.newnew1} \e
and
  \ee \|\w_h\|_{H^1_{\ast}(\o)}+h^{-1}\|{\bf R}_h\|_{L^2_{\omega}(\o)}\stlm
\log^2(1/h)\rho^{m_0}(h)\|\v_h\|_{H^{\ast}(\c,~\o)},\label{stab3.newnew2}
  \e
  where the constant $m_0$ is independent of both $h$ and the coefficients $\alpha$ and $\beta$.
\end{theorem}

The following result can be obtained directly from Theorem \ref{thm:main}.
\begin{corollary} Assume that the coefficients satisfy $c_0\al\leq\be\leq C_0\al$ for two given positive constants $c_0$ and $C_0$. Then,  every function $\v_h\in V_h^{\ast}(\Omega)$ admits
the decomposition (\ref{newnew3.01}), which satisfies the stability estimates (\ref{stab3.newnew1}) and (\ref{stab3.newnew2}) with $\mu=\max\{C_0,1/c_0\}$.
\end{corollary}

It is natural to wonder whether additional constraints are indeed necessary for a function $\v_h$ to admit a stable
regular decomposition for the case with {\it thorny vertices}. The following conclusion can be obtained by considering a counterexample.\\
{\bf Proposition 3.3} The space $V^{\ast}_h(\Omega)$ in Theorem \ref{thm:main} cannot be replaced by $V_h(\Omega)$ itself.\\
\no{\it Proof}. We need only to construct a counterexample. To this end, we consider the ``checkerboard" domain (see Figure 1). Let $\varepsilon\ll 1$ be a very small positive number.
We define the coefficients $\al$ and $\be$ as follows: $\al(\x)=1$ and $\be(\x)=\varepsilon$ on $\Omega_1\cup\Omega_2$,
and $\al(\x)=\be(\x)=\varepsilon$ on $\Omega_3\cup\Omega_4$. For this example, there is only one
{\it thorny vertex} $\vv$ at the center of $\Omega$ and $\Im_{\vv}^{\ast}$ contains two cubes $\Omega_1$ and $\Omega_2$, which implies that $n_{\ast}=1$.
It is easy to see that the coefficients in this example satisfy {\bf Condition A}, {\bf Condition B} and {\bf Condition C} in the definition of the class $\mathbb{P}(\mu)$.

For $i=1,2$, let $\phi_{h,i}\in Z_h(\Omega_i)$ be
a nodal finite element function satisfying $\phi_{h,i}=0$ on $\partial \Omega_i\cap\partial\Omega$. Construct an edge finite element function $\v_h$ such that:
(i) $\v_h=\nabla\phi_{h,i}$ on $\bar{\Omega}_i$ ($i=1,2$); (ii) $\lambda_e(\v_h)=0$ for any edge $e$ in $\Omega_3\cup\Omega_4$; and (iii) $\v_h$ has continuous tangential components on
every fine edge on $(\partial\Omega_3\cup\partial\Omega_4)\cap(\partial \Omega_1\cup\partial \Omega_2)$. It is easy to see that $\v_h\in V_h(\Omega)$
even if $\phi_{h,1}\not=\phi_{h,2}$ at $\vv$ (an edge finite element function may be discontinuous at a node). Note that $\v_h$ does not vanish on $\Omega_3\cup\Omega_4$
because $\v_h$ has non-zero degrees of freedom on $(\partial\Omega_3\cup\partial\Omega_4)\cap(\partial \Omega_1\cup\partial \Omega_2)$.
In the following, we explain that the function $\v_h$ must satisfy an additional constrain if it admits a regular decomposition
satisfying all the requirements in Theorem \ref{thm:main} (with $\rho(h)=1$).

We assume that $\v_h$ admits the regular decomposition (\ref{newnew3.01}). For convenience, set $G=\Omega_1\cup\Omega_2$. By the definition of $\v_h$, we have $\c\,\v_h=\0$ on $G$.
Then, the estimate (\ref{stab3.newnew2}) implies that
$$ |\w_h|_{1,G}+h^{-1}\|{\bf R}_h\|_{0,G}\stlm\log(1/h)\rho^m(h)\varepsilon^{\12}(\|\v_h\|_{0,G}+\|\v_h\|_{H^{\ast}(\c,\Omega_3\cup\Omega_4)}). $$
For a fixed $h$, let $\varepsilon\rightarrow 0^+$, the above inequality infers that $\w_h={\bf R}_h=\0$ on $G$
(note that $\w_h$ vanishes on $\partial\Omega\cap\partial G$). Therefore, by the regular decomposition, we have $\v_h=\nabla p_h$ on $G$ with $p_h\in Z_h(\Omega)$.
For $i=1,2$, set $p_{h,i}=p_h|_{\Omega_i}$ and let $\ff_i\subset\partial \Omega_i$ be a face containing $\vv$ as one of its vertex. For this example, we have $\ff_i\cap\partial\Omega\not=\emptyset$, so we
can choose a vertex $\vv_i$ on $\partial\ff_i\cap\partial\Omega$ such that $p_{h,i}$ vanishes at $\vv_i$. In the following, we calculate the arc-length integrals of $\v_h\cdot\t_{\partial\ff_i}$ on $\partial\ff_i$.

Without loss of generality, we assume that the arc-length coordinate of the point $\vv_i$ is just $0$. For convenience, we use $t_{\vv}$ to denote the arc-length coordinate of the center vertex $\vv$.
By the condition $\v_h=\nabla p_h$ on $\Omega_i$, we deduce that
$$ \int_0^{t_{\vv}}\v_h\cdot\t_{\partial\ff_i}ds=\int_0^{t_{\vv}}\nabla p_{h,i}\cdot\t_{\partial\ff_i}ds=p_{h,i}(t_{\vv})~~~~(i=1,2).$$
Since $p_h\in Z_h(\Omega)$, we have $p_{h,1}(t_{\vv})=p_{h,2}(t_{\vv})$. Thus
$$  \int_0^{t_{\vv}}\v_h\cdot\t_{\partial\ff_1}ds=\int_0^{t_{\vv}}\v_h\cdot\t_{\partial\ff_2}ds. $$
Namely, the function $\v_h$ must satisfy the condition ${\mathcal F}\v_h=0$ with
\ee {\mathcal F}\v_h=\int_0^{t_{\vv}}\v_h\cdot\t_{\partial\ff_1}ds-\int_0^{t_{\vv}}\v_h\cdot\t_{\partial\ff_2}ds. \label{example-constraint}
\e
Then, this proposition is proved.
\hfill $\Box$

{\bf Proposition 3.3} tells us that, for the case with {\it thorny vertices}, additional constraints are indeed necessary for a function $\v_h$ that admits a stable regular decomposition.
For this particular example, the functional in the subspace $V_h^{\ast}(\Omega)$ can be replaced by ${\mathcal F}$ given in (\ref{example-constraint}), but the conclusion
is not valid for general situations. It is easy to see that the condition in $V_h^{\ast}(\Omega)$ is also satisfied for this particular $\v_h$ constructed in the above example.

\subsection{The case without thorny vertex}

In this subsection, we consider a particular case without {\it thorny vertex}, which has some connection with the {\it quasi-monotonicity assumption} introduced in \cite{Pet2002}.

The original definition of the quasi-monotonicity assumption is somewhat complicated (see Definitions 4.1 and 4.6
in Subsection 4.1 of \cite{Pet2002}). Thanks to the introduction of the {\it upper intersection set}, the {\it quasi-monotonicity assumption} can be described
in a succinct manner: the distribution of the coefficients $\al_1,\al_2,\cdots,\al_{N_0}$ satisfies the {\it quasi-monotonicity assumption} if,
for each subdomain $\Omega_k$, the {\it upper intersection set} ${\mathcal S}_k$ is the empty set or a union of some faces of $\Omega_k$, i.e., there is no isolated edge or
isolated vertex in ${\mathcal S}_k$.

In order to compare the considered case with the {\it quasi-monotonicity assumption}, we would like to introduce a new concept.
\begin{definition}\label{definition2} We say that the distribution of the coefficients $\al_1,\al_2,\cdots,\al_{N_0}$ satisfies the {\it generalized quasi-monotonicity assumption} if,
for each subdomain $\Omega_k$, the {\it upper intersection set} ${\mathcal S}_k$ is the empty set or a union of some faces and edges of $\o_k$, i.e., there is no isolated vertex in ${\mathcal S}_k$.
\end{definition}
\vskip 2mm

It is clear that the {\it generalized quasi-monotonicity assumption} is weaker than the {\it quasi-monotonicity assumption}.

The following result comes from {\bf Proposition 3.2} and Definition \ref{definition2}.\\
{\bf Proposition 3.4}. There is no {\it thorny vertex} if and only if the distribution of the coefficients $\al_1,\al_2,\cdots,\al_{N_0}$ satisfies {\it generalized quasi-monotonicity assumption}.

This proposition implies that, when the coefficients $\al_1,\al_2,\cdots,\al_{N_0}$ satisfies {\it generalized quasi-monotonicity assumption}, we have $\Im^{\ast}=\emptyset$ and $\Im^c_{\ast}=\{\Omega_k\}_{k=1}^{N_0}$.

\begin{definition}\label{classM} Let $\mu\geq 1$ be a given constant. We say that the coefficients $\{\al_k\}$ and $\{\be_k\}$ belong to the problem class $\mathbb{P}_M(\gamma)$ if the distribution of the coefficients
$\al_1,\al_2,\cdots,\al_{N_0}$ satisfies the {\it generalized quasi-monotonicity assumption}, and the coefficients $\{\al_k\}$ and $\{\be_k\}$ satisfy {\bf Condition A} in Definition \ref{class} and
the following condition:

{\bf Condition B$'$}. For every $\Omega_k$, on which the coefficients do not satisfy the relation $\gamma^{-1}\al_k\leq \be_k\leq\gamma\al_k$, the {\it upper intersection set} ${\mathcal S}_k$
(if it is non-empty) meets one of the following two requirements:

(i) if $\be_k>\mu\al_k$, then ${\mathcal S}_k$ is a (possibly non-connected) union of some {\it faces} of $\Omega_k$;

(ii) if $\be_k<\mu^{-1}\al_k$, then ${\mathcal S}_k$ is a {\it connected} union of some faces and edges of $\Omega_k$.\\

\end{definition}


For the case considered in this subsection, we have the following result
\begin{theorem}\l{thm:main1} Assume that the coefficients $\{\al_k\}$ and $\{\be_k\}$ belong to the problem class $\mathbb{P}_M(\gamma)$.
Then, any function $\v_h\in V_h(\o)$ admits a decomposition of the form \ee
  \v_h=\nabla p_h+\r_h\w_h+{\bf R}_h \label{newnew3.001}\e for some
  $p_h\in Z_h(\o)$ and $\w_h\in (Z_h(\o))^3$ and ${\bf R}_h\in
  V_h(\o)$. Moreover, $p_h$, $\w_h$ and ${\bf R}_h$ satisfy the estimates
  \ee \|p_h\|_{H^1_{\beta}(\o)}\stlm\rho^{m_0}(h)\|\v_h\|_{H^{\ast}(\c,~\o)}\label{stab3.newnew01} \e
  and\ee \|\w_h\|_{H^1_{\ast}(\o)}+h^{-1}\|{\bf R}_h\|_{L^2_{\omega}(\o)}\stlm
  \rho^{m_0}(h)\|\v_h\|_{H^{\ast}(\c,~\o)}, \l{newnew3.02} \e
  where the constant $m_0$ ($\geq 2$) is independent of both $h$ and the coefficients $\omega$ and $\beta$.
\end{theorem}

\begin{remark} If the {\it upper intersection set} ${\mathcal S}_k$ is a union of connected Lipschitz unions of some faces of $\Omega_k$ for each subdomain $\Omega_k$,
then the factor $\rho^{m_0}(h)$ in (\ref{stab3.newnew01}) and (\ref{newnew3.02}) can be dropped.
\end{remark}

The following result can be obtained directly from Theorem \ref{thm:main1}.
\begin{corollary} Assume that the coefficients satisfy the {\it generalized quasi-monotonicity assumption} and the relation
$c_0\al\leq\be\leq C_0\al$ for two given positive constants $c_0$ and $C_0$. Then, every function $\v_h\in V_h(\Omega)$ admits
the decomposition (\ref{newnew3.001}), which satisfies the stability estimates (\ref{stab3.newnew01}) and (\ref{newnew3.02}) with $\mu=\max\{C_0,1/c_0\}$.
\end{corollary}

The assumptions in Theorem \ref{thm:main1} are satisfied in relevant settings as in the example constructed in Subsection 6.1 of the paper \cite{Kol2009}.
This example reads as follows.

\begin{example}\label{example2} Let $\Omega = [0,1]^{3}$ be the unit cube, and be divided into two
polyhedra $D_1$ and $D_2$ (see Figure 5). There are two choices of
the coefficients: (a) $\alpha\equiv1$ on $\o$, $\beta=1$ on $D_1$
and $\beta= 10^{k}$ (with $k = -8,\cdots, 8$) on $D_2$; (b)
$\beta\equiv1$ on $\o$, $\alpha=1$ on $D_1$ and $\alpha=10^{k}$
(with $k = -8, \cdots, 8$) on $D_2$. A uniform triangulation is used
on $\o$.
\end{example}

\begin{center}
\includegraphics[width=5cm,height=5cm]{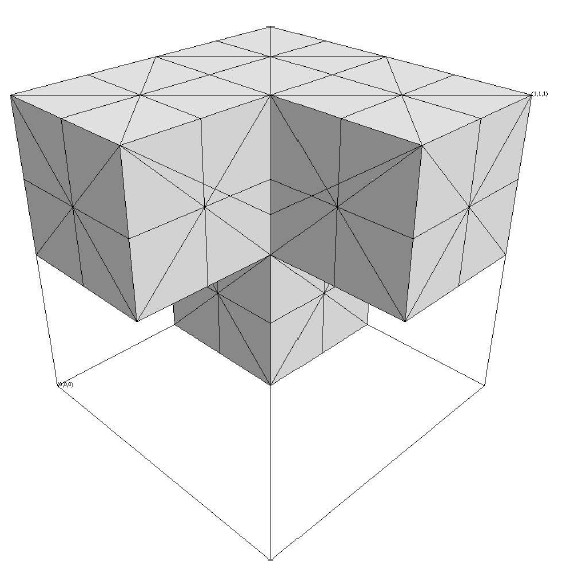}\quad
\includegraphics[width=5cm,height=5cm]{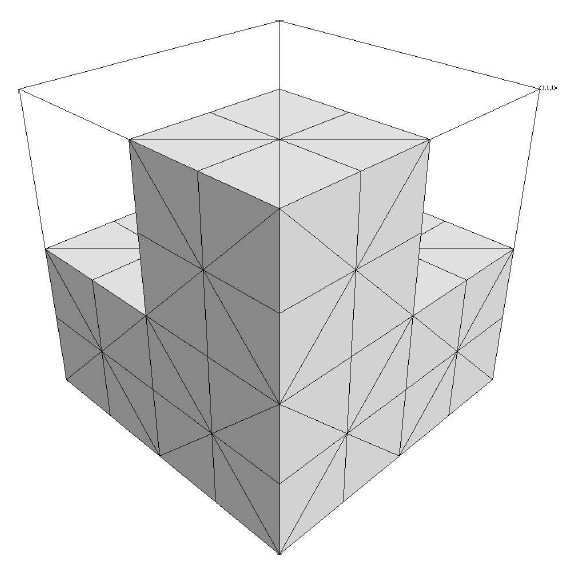}

\centerline{Figure 5: The unit cube split into two symmetrical
regions (left: $D_1$; right: $D_2$)}
\end{center}
\vskip 2mm

For this example, regarding $D_k$ as $\Omega_k$, the {\it upper intersection set} ${\mathcal S}_k$ is a connected union of faces for $k=1,2$,
which means that the {\it generalized quasi-monotonicity assumption} is satisfied.
For Case (a), we have $\al_1=\al_2$ and $\be_1=\al_2<\beta_2$ when $k=1,\cdots,8$; we have $\al_2=\al_1$ and $\be_2<\beta_1=\al_1$ when $k=-8,\cdots,0$, so {\bf Condition A}
in Definition \ref{class} and {\bf Condition B$'$} in Definition \ref{classM} are met for $\mu=1$.
For Case (b), when $k=-8,\cdots,0$ we have $\al_2\leq \al_1$ and $\be_2=\al_1=\be_1$; when $k=1,\cdots,8$ we have $\al_1<\al_2$ and $\be_1=\al_2=\be_2$, so
{\bf Condition A} in Definition \ref{class} and {\bf Condition B$'$} in Definition \ref{classM} are also satisfied for $\mu=1$.
Then the coefficients in this example belong to the problem class $\mathbb{P}_M(\gamma)$ with $\mu=1$. Therefore, Theorem \ref{thm:main1} holds and the resulting decomposition is uniformly stable.


\section{Application to analysis of the
HX preconditioner with jump coefficients}\l{sec:appl}
\setcounter{equation}{0}

In this section, we apply the discrete weighted regular
decompositions described in Theorems~\ref{thm:main} and \ref{thm:main1} to analyze the convergence
of the HX preconditioner for the case with jump coefficients.

\subsection{The HX preconditioner}
\l{sec:HX preconditioner}

First, we recall the HX preconditioner proposed in \cite{HX} for the linear system of equations arising
from the finite element Galerkin discretization of (\ref{1.1}).

The discrete variational problem of (\ref{1.1}) is: to find $\u_h\in V_h(\o)$ such that
\ee
(\al\c\,\u_h,\c\,\v_h)+(\beta\u_h,\v_h)=(\f,\v_h),~~~\forall\v_h\in
V_h(\o).\l{6.1} \e As usual, we can rewrite this in the operator form
\ee {\bf A}_h\u_h=\f_h,\l{6.02}\e with ${\bf A}_h:V_h(\o)\rightarrow
V_h(\o)$ defined by
$$ ({\bf A}_h\u_h,\v_h):=(\al\c\,\u_h,\c\,\v_h)+(\beta\u_h,\v_h),~~\u_h,\v_h\in V_h(\o).$$

Let ${\bf \Delta}_h: (Z_h(\o))^3\rightarrow (Z_h(\o))^3$ be the
discrete elliptic operator defined by
$$ (-{\bf \Delta}_h\v,\w):=(\al\nabla\v,\nabla\w)+(\be\v,\w),~~\v,\w\in (Z_h(\o))^3,$$
and let ${\bf A}^{\nabla}_h:\nabla(Z_h(\o))\rightarrow \nabla(Z_h(\o))$
be the restriction of ${\bf A}_h$ on the space $\nabla(Z_h(\o))$,
whose action can be implemented by solving a Laplace equation.
Additionally, let ${\bf J}_h:V_h(\o)\rightarrow V_h(\o)$ denote the standard
Jacobi smoother of ${\bf A}_h$, where ${\bf J}_h$ corresponds to the diagonal preconditioner
for the stiffness matrix of ${\bf A}_h$ and was characterized by (7.1) of \cite{HX}.. Then the HX preconditioner ${\bf B}_h$ of
${\bf A}_h$ can be defined by
$$ {\bf B}_h:={\bf J}_h^{-1}+\hat{\bf r}_h(-{\bf \Delta}_h)^{-1}\hat{\bf
r}_h^{\ast}+{\bf T}_h({\bf A}^{\nabla}_h)^{-1}{\bf T}_h^{\ast}, $$ where
 $\hat{\bf r}_h: (Z_h(\o))^3\rightarrow V_h(\o)$ is
the restriction of the interpolation operator $\r_h$ on
$(Z_h(\o))^3$, and ${\bf T}_h^{\ast}: V_h(\o)\rightarrow
\nabla(Z_h(\o))$ is the $L^2$ projector.

It has been shown that the condition number of the preconditioned system is independent of the meshwidth $h$ (see Theorem 7.1 in Subsection 7.1 of \cite{HX}):
$$ cond({\bf B}_h{\bf A}_h)\leq C $$
with $C>0$ independent of $h$. However, it is unclear how the constant $C$ depends on the jumps in the
coefficients $\al$ and $\be$.

\subsection{Convergence of the HX preconditioner for the case with
  jump coefficients}

In this subsection, we derive a new convergence result for the
preconditioner ${\bf B}_h$ in the case where the
coefficients $\al$ and $\be$ have large jumps. This new result follows from Theorems \ref{thm:main} and \ref{thm:main1}.

Let $n_{\ast}$ and $V_h^{\ast}(\Omega)$ be as defined in Subsection 3.2. We use $\lambda_{n_{\ast}+1}({\bf B}_h^{-1}{\bf A}_h)$ to denote the minimal eigenvalue of
the restriction of ${\bf B}_h^{-1}{\bf A}_h$ on the subspace $V^{\ast}_h(\Omega)$, and define $\kappa_{n_{\ast}+1}({\bf B}_h^{-1}{\bf A}_h)$
as the {\it reduced condition number} (see \cite{XuZ2008}) of ${\bf B}^{-1}_h{\bf A}_h$ associated with the subspace $V^{\ast}_h(\Omega)$. Namely,
$$\kappa_{n_{\ast}+1}({\bf B}_h^{-1}{\bf A}_h):={\lambda_{\max}({\bf B}_h^{-1}{\bf A}_h)\over \lambda_{n_{\ast}+1}
({\bf B}_h^{-1}{\bf A}_h)}. $$
Under the framework introduced in \cite{XuZ2008}, the convergence rate of the PCG method with the preconditioner ${\bf B}_h$ for solving system (\ref{6.02})
is determined by the {\it reduced condition number} $\kappa_{n_{\ast}+1}({\bf B}_h^{-1}{\bf A}_h)$, and the number of iteration required to achieve a given accuracy of the approximation
weakly depends on the codimension $n_{\ast}$, instead of the space $V^{\ast}_h(\o)$ itself (i.e., the choice of the functionals $\{{\mathcal F}_{i,r}\}$).
If there is no {\it thorny vertex}, then
we have $n_{\ast}=0$ and $V_h^{\ast}(\Omega)=V_h(\Omega)$, which means that $\kappa_{n_{\ast}+1}({\bf B}_h^{-1}{\bf A}_h)$ is just the standard condition number $cond({\bf B}_h^{-1}{\bf A}_h)$.
In this part, we focus on the estimate of $\kappa_{n_{\ast}+1}({\bf B}_h^{-1}{\bf A}_h)$.

Let the function $h\mapsto \rho(h)$ be defined as in Subsection 3.2, so it takes a value of $1$ or $\log(1/h)$. The following results can be proved in the standard manner
(the details can be found in Subsection 4.2 of \cite{Hu2-preprint}).
\begin{theorem}\label{HX} Assume that the coefficients $\{\al_k\}$ and $\{\be_k\}$ belong to the problem class $\mathbb{P}(\gamma)$.
Then, the {\it reduced condition number} of the preconditioned system satisfies
\ee \kappa_{n_{\ast}+1}({\bf B}^{-1}_h{\bf A}_h)\stlm\log^2(1/h)\rho^{m_0}(h).\l{6.2} \e
When the coefficients $\{\al_k\}$ and $\{\be_k\}$ belong to the problem class $\mathbb{P}_M(\gamma)$ (there is no {\it thorny vertex}), we have
\ee
cond({\bf B}^{-1}_h{\bf A}_h)\stlm\rho^{m_0}(h).\l{6.2newnew} \e
In particular, if every set ${\mathcal S}_k$ is a union of connected Lipschitz unions of some faces of $\Omega_k$ for each subdomain $\Omega_k$,
we have
\ee
cond({\bf B}^{-1}_h{\bf A}_h)\stlm 1.\l{6.2new} \e
\end{theorem}
\no{\it Proof}. We need only to prove the inequality (\ref{6.2}) by Theorem \ref{thm:main}, and we can derive (\ref{6.2newnew}) and (\ref{6.2new}) in the same manner
(but using Theorem \ref{thm:main1}). For any $\v_h\in V^{\ast}_h(\o)$, let $\w_h\in (Z_h(\o))^3$, $p_h\in Z_h(\o)$ and $\R_h\in V_h(\o)$ be defined by the
decomposition in Theorem \ref{thm:main}, i.e.,
\ee \v_h=\nabla p_h+\r_h\w_h+\R_h. \l{6.5}\e By the auxiliary space
technique for the construction of preconditioner (refer to
\cite{HX}), it suffices to verify that
\ee (\Delta_h(\nabla p_h),\nabla p_h)+({\bf
\triangle}_h\w_h,\w_h)+({\bf J}\R_h,\R_h)\stlm\log^2(1/h)\rho^{m_0}(h)({\bf A}_h\v_h,\v_h).\l{6.3} \e

From (\ref{stab3.newnew1}), we have \beq
(\Delta_h(\nabla p_h),\nabla p_h)&\stl&\|\be^{1\over 2}\nabla
p_h\|^2_{0,\o}\stlm\log^2(1/h)\rho^{m_0}(h)(\|\al^{1\over
2}\c\,\v_h\|^2_{0,\o}+\|\be^{1\over 2}\v_h\|^2_{0,\o})\cr
&\stlm&\log^2(1/h)\rho^{m_0}(h)({\bf A}_h\v_h,\v_h).\l{6.8}\eq
It follows from (\ref{stab3.newnew2}) that  \beq ({\bf
\triangle}_h\w_h,\w_h)&=&\|\w_h\|^2_{H^1_{\ast}(\o)}\stlm\log^2(1/h)\rho^{m_0}(h)\|\v_h\|^2_{H^{\ast}(\c;~\o)}\cr
&\stlm&\log^2(1/h)\rho^{m_0}(h)({\bf A}_h\v_h,\v_h)\l{6.7}\eq
and
 \beq ({\bf
J}{\bf R}_h,{\bf R}_h&\stl& h^{-2}\|\al^{1\over
2}{\bf R}_h\|^2_{0,\o}\stlm \log^2(1/h)\rho^{m_0}(h)\|\al^{1\over
2}\c\,\v_h\|^2_{0,\o}\cr&\stlm&\log^2(1/h)\rho^{m_0}(h)({\bf
A}_h\v_h,\v_h).\l{6.9}\eq Then, the inequality (\ref{6.3}) is a direct consequence
of (\ref{6.8})-(\ref{6.9}).
\hfill$\Box$
\begin{remark} In most applications, there are relatively few {\it thorny vertices}, even if many different media are involved in the entire physical domain.
Thus, the codimension $n_{\ast}$ is often a small positive integer independent of both $h$ and the jumps in the coefficients $\al$ and $\be$ (but it depends on the
distribution of the coefficient $\{\al_k\}$). Then, by Theorem \ref{HX} and the framework introduced in \cite{XuZ2008},
the PCG method with the preconditioner ${\bf B}_h$ for solving the system (\ref{6.02}) achieves fast convergence, provided that the coefficients $\{\al_k\}$ and $\{\be_k\}$
belong to the class $\mathbb{P}(\gamma)$ with a mild parameter $\mu$. In particular, the convergence rate is nearly optimal
up to a poly-logarithmic factor in the meshwidth $h$ and is not destroyed by large jumps of the coefficients $\al$ and $\be$ across the interface between any two neighboring subdomains.
\end{remark}

\section{Proofs} In this section, we prove the main results Theorems~\ref{thm:main} and ~\ref{thm:main1}.
The basic idea is to transform the general case involved in Theorems~\ref{thm:main} into the particular case considered in Theorems~\ref{thm:main1}, so we first
prove Theorem \ref{thm:main1}.

\subsection{Analysis for the case satisfying the generalized quasi-monotonicity assumption}
\setcounter{equation}{0}
In this subsection, we focus on the proof of Theorem \ref{thm:main1}.
Throughout this subsection, we always assume that the coefficients $\{\al_k\}$ and $\{\be_k\}$ belong to the problem class $\mathbb{P}_M(\gamma)$.

From now on, we say that two subdomains $\o_r$ and $\o_{r'}$ do not intersect if $\bar \o_r\cap \bar\o_{r'}=\emptyset$; otherwise, the two subdomains intersect each other.
For a polyhedron $G$ in $\{\Omega_k\}$, we use $\Xi_G$ to denote the union of all the polyhedra that belong to $\{\Omega_k\}$ and intersect with $G$. The following result will be used
repeatedly.
\begin{lemma}\label{extension-new} Let $G$ be a polyhedron in $\{\Omega_k\}$ and $\Gamma$ be the union of some faces and edges of $G$. Then, there exists an extension $E_h$ mapping $Z_h(G)$ into $Z_h(\Omega)$ such that,
for any function $\phi_h\in Z_h(G)$ that vanishes on $\Gamma$, the function $E_h\phi_h$ satisfies the following conditions: (1) $E_h\phi_h=\phi_h$ on $\bar{G}$; (2) $supp~E_h\phi_h\subset \Xi_G$; (3) when $G'\subset\Xi_G$ and
$\bar{G}'\cap \bar{G}\subset\Gamma$, we have $E_h\phi_h=0$ on $G'$; (4) it fulfills the stability estimates
\ee
\|E_h\phi_h\|_{1,\Omega}\leq C\log (1/h)\|\phi_h\|_{1,G}\quad\mbox{and}\quad\|E_h\phi_h\|_{0,\Omega}\leq C\|\phi_h\|_{0,G}.\label{extension:4.0001}
\e
In particular, if $\Gamma$ does not contain any edge, the factor $\log (1/h)$ in the above inequality can be removed.
\end{lemma}

\no{\it Proof}. When $\Gamma$ is a union of only some faces of $G$, the results can be proved as in the proof of Lemma 4.5 in \cite{Klaw2008}.

For the case that $\Gamma$ contains edges, we define the extension in the same manner as $\tilde{\w}_{h,1}$ in the proof of Theorem 4.1 in \cite{Hu1-2017}.
Assume that $G$ has $n_{f}$ faces, which are denoted by $\ff_1,\cdots,\ff_{n_f}$.
Set $\ff^{\partial}=\cup_{j=1}^{n_f}\partial\ff_j$. For each $\ff_j$, let $\vartheta_{\ff_j}$ be the finite element
function defined in \cite{Cas1996} and \cite{Wid05}. This function satisfies
$\vartheta_{\ff_j}(\x)=1$ for each node $\x\in\bar{\ff}_j\backslash\partial\ff_j$,
$\vartheta_{\ff_j}(\x)=0$ for $\x\in\partial G\backslash\ff_j$ and $0\leq\vartheta_{\ff_j}\leq 1$ on
$G$. Let $\pi_h$ denote the standard interpolation operator into
$Z_h(G)$, and define $\phi^{\ff_j}_h=\pi_h(\vartheta_{\ff_j}\phi_h)$ ($j=1,\cdots,n_f$).

When $\ff_j\subset\Gamma$, we define the extension $\tilde{\phi}^{\ff_j}_h$ of $\phi^{\ff_j}_h$ as the natural zero extension because $\phi^{\ff_j}_h=0$ on $\Gamma$.
We need only to consider the faces $\ff_j$ that are not contained in $\Gamma$. Let $G_j\subset\Xi_G$ be the polyhedron having the common face $\ff_j$ with $G$. As in
Lemma 4.5 of \cite{Klaw2008}, we can show there exists
an extension $\tilde{\phi}^{\ff_j}_h$ of $\phi^{\ff_j}_h$ such that $\tilde{\phi}^{\ff_j}_h\in Z_h(\Omega)$; $\tilde{\phi}^{\ff_j}_h=\phi^{\ff_j}_h$ on $\bar{G}$; $\tilde{\phi}^{\ff_j}_h$ vanishes
on $\Omega\backslash(G\cup\ff_j\cup G_j)$; $\tilde{\phi}^{\ff_j}_h$ is stable with both $H^1$ norm and $L^2$ norm. Define
$$ E_h\phi_h=\phi^{\partial}_h+\sum_{j=1}^{n_f}\tilde{\phi}^{\ff_j}_h, $$
where $\phi^{\partial}_h\in Z_h(\Omega)$ denotes the zero extension of the restriction of $\phi_h$ on $\ff^{\partial}$. Then, the extension $E_h\phi_h$ meet all the requirements in this lemma.
\hfill $\Box$

In the rest of this section, for a nodal finite element function $\phi_h$ we always
use $\tilde{\phi}_h$ to denote its extension as defined by Lemma \ref{extension-new}. For convenience, such an extension is simply called
a {\it stable} extension.

For a subdomain $\Omega_k$, let ${\mathcal S}_k$ be the {\it upper intersection set} defined by (\ref{intersection}).
In the analysis below, we will consider four cases simultaneously (where $\mu\geq 1$ denotes a given constant in the definition of the class $\mathbb{P}_M(\gamma)$):

Case (i): $\be_k<\mu^{-1}\al_k$ and ${\mathcal S}_k$ is the connected union of faces of $\Omega_k$;

Case (ii): $\be_k<\mu^{-1}\al_k$ and ${\mathcal S}_k$ is the connected union of faces and edges of $\Omega_k$;

Case (iii): $\al_k<\mu^{-1}\be_k$ and ${\mathcal S}_k$ is the union (possibly non-connected) of some faces of $\Omega_k$;

Case (iv): $\mu^{-1}\beta\leq\al_k\leq\mu\be_k$ and ${\mathcal S}_k$ is the union (possibly non-connected) of some faces and edges of $\Omega_k$.

For Cases (i) and (ii), we use the first particular conclusion of Theorem 3.1 in Section 3 of \cite{Hu1-2017}
to build a regular decomposition on $\Omega_k$. For Cases (iii) and (iv), we respectively use the second particular conclusion and the general conclusion of Theorem 3.1 in Section 3
of \cite{Hu1-2017} to do so.

For convenience, for $\v_h\in V_h(\Omega)$, we define
two general ``norms", $\|\v_h\|_{\ast,\Omega_k}$ and $\|\v_h\|_{\#,\Omega_k}$, as follows: $\|\v_h\|_{\ast,\Omega_k}=\|\c~\v_h\|_{0,\Omega_k}$ for Cases (i) and (ii), or
$\|\v_h\|_{\ast,\Omega_k}=\|\v_h\|_{\c,\Omega_k}$ for Cases (iii) and (iv); $\|\v_h\|_{\#,\Omega_k}=\|\v_h\|_{0,\Omega_k}$ for Case (iii), or
$\|\v_h\|_{\#,\Omega_k}=\|\v_h\|_{\c,\Omega_k}$ for the other cases.

As in Subsection 3.2, we define the positive function $\rho_k(h)$ as follows: $\rho_k(h)=1$ when
${\mathcal S}_k$ is an empty set or is a {\it union of connected Lipschitz unions of some faces} of $\Omega_k$
; otherwise, $\rho_k(h)=\log(1/h)$ (if ${\mathcal S}_k$ contains an isolated edge or a non-Lipschitz union of faces).

\subsubsection{A decomposition for edge element functions}\l{sec:decomp}
\setcounter{equation}{0}
In this part, we build a suitable regular decomposition for functions $\v_h\in V_h(\Omega)$.


The basic ideas come from \cite{HuShuZou2013} and can be described as follows.
We first decompose $\{\o_r\}_{r=1}^{N_0}$ into a union of non-empty subsets $\Sigma_1$, $\cdots$, $\Sigma_m$ satisfying the following conditions:
(1) any two polyhedra in the same subset $\Sigma_l$ do not intersect each other; (2) for any two polyhedra $\Omega_{r_l}$ and $\Omega_{r_j}$ belonging
to two different subsets $\Sigma_l$ and $\Sigma_{j}$, respectively, with $l<j$, we have $\al_{r_l}\geq\al_{r_j}$ if $\Omega_{r_l}$ and $\Omega_{r_j}$ intersect each other.
Then we in turn construct the desired decomposition (\ref{newnew3.001}) on the subdomains in the subsets $\Sigma_1$, $\cdots$, $\Sigma_m$.

Specific definitions of $\{\Sigma_l\}$ were described in Section 5 of \cite{HuShuZou2013}, from which we see that the value of the positive integer $m$ only depends upon the distribution
of the coefficients $\{\al_k\}$. Without loss of generality, we assume that
$$ \Sigma_l=\{\o_{n_{l-1}+1},~\o_{n_{l-1}+2},\cdots,\o_{n_l}\} $$
with $n_0=0$ and $n_l>n_{l-1}$ $(l=1,\cdots,m)$. It is clear that $\Sigma_l$ contains $(n_l-n_{l-1})$ polyhedra.

\smallskip

We are now ready to construct the desired decomposition (\ref{newnew3.001}) for any $\v_h$ in $V_h(\Omega)$. We do so by three steps absorbing the idea of the inductive approach.
In the following, we always use $\v_{h,r}$ to denote the restriction of $\v_h$ on a subdomain $\Omega_r$.

\smallskip
{\bf Step 1}: Decompose $\v_h$ on all the polyhedra in $\Sigma_1$.

For $\o_r\in\Sigma_1$ (i.e., $1\leq r\leq n_1$), the function $\v_{h,r}$ vanishes on the set $\hat{\mathcal S}_r=\partial\o_r\cap\partial\o$,
which is a union of some faces of $\o_r$.
Using Theorem 3.1 in Section 3 of \cite{Hu1-2017} (regarding $\v_{h,r}$ and $\hat{\mathcal S}_r$ as $\v_h$ and $\Gamma$, respectively), we have the decomposition:
\ee
\v_{h,r}=\nabla p_{h,r}+\r_h\w_{h,r}+{\bf R}_{h,r}\,,\l{4.001} \e
where $p_{h,r}\in Z_h(\o_r)$, $\w_{h,r}\in (Z_h(\o_r))^3$ and ${\bf
R}_{h,r}\in V_h(\o_r)$, and they vanish on $\hat{\mathcal S}_r$. Moreover, we have
\ee
\|\w_{h,r}\|_{1,\o_r}+h^{-1}\|{\bf R}_{h,r}\|_{0,\o_r}\stl\|\v_{h,r}\|_{\ast,\o_r}\label{4.009}
\e
and
\ee
\|\w_{h,r}\|_{0,\o_r}+\|p_{h,r}\|_{1,\o_r}\stl\|\v_{h,r}\|_{\#,\o_r}.\l{stab:5.05newnew} \e

Let $\tilde{p}_{h,r}\in Z_h(\o)$ and $\tilde{\w}_{h,r}\in (Z_h(\o))^3$ be the {\it stable} extensions of $p_{h,r}$ and $\w_{h,r}$,
respectively. Moreover, let $\tilde{{\bf R}}_{h,r}\in
V_h(\o)$ denote the natural zero extensions of ${\bf R}_{h,r}$. Then, we set
\ee
\tilde{\v}_{h,r}=\nabla
\tilde{p}_{h,r}+\r_h\tilde{\w}_{h,r}+\tilde{{\bf R}}_{h,r}~~ \m{for
all $r$ such that} ~\o_r\in\Sigma_1.\l{4.014} \e

{\bf Step 2}: Decompose $\v_h$ on all the polyhedra in $\Sigma_2$.

Consider a subdomain $\o_r$ from $\Sigma_2$. For ease of notation,  we introduce the index set
$$ \Lambda_r^1=\{~i~; ~1\leq i\leq n_1 \m{~such that ~}\partial\o_i\cap\partial\o_r\not=\emptyset\},\quad\Omega_r\in\Sigma_l~~(l\geq 2).
$$
Define the remainder
\ee
\v^{*}_{h,r}=\v_{h,r}-
\sum\limits_{i\in \Lambda_r^1}\tilde{\v}_{h,i}~~\mbox{on}~\o_r \l{decom5.new1} \e
and
$$ \hat{\mathcal S}_r=\bigcup_{i\in \Lambda_r^1}(\partial\Omega_r\cap\partial\Omega_i)\bigcup(\partial\Omega_r\cap\partial\Omega).  $$

From (\ref{4.001}), (\ref{4.014}) and (\ref{decom5.new1}), we have $\la_e(\v^{*}_{h,r})=0$ for $e\subset\hat{\mathcal S}_r$. It is clear that $\hat{\mathcal S}_r$
is a subset of the {\it upper intersection set} ${\mathcal S}_r$. Thus, by {\bf Condition B$'$} in the definition of the class $\mathbb{P}_M(\mu)$, the set $\hat{\mathcal S}_r$
is a union of faces and edges. Then, using Theorem 3.1 in Section 3 of \cite{Hu1-2017} again, there exist $p^{*}_{h,r}\in Z_h(\o_r)$,
$\w^{*}_{h,r}\in (Z_h(\o_r))^3$ and ${\bf R}^{*}_{h,r}\in V_h(\o_r)$ such that they have zero degrees of freedom on $\hat{\mathcal S}_r$ and satisfy
\ee
\v^{*}_{h,r}=\nabla p^{*}_{h,r}+\r_h\w^{*}_{h,r}+{\bf
R}^{*}_{h,r}~~ \mbox{on}~~\o_r\,.\l{4.010}
\e
Moreover, we have (according to specific case)
\ee
\|\w^{*}_{h,r}\|_{1,\o_r}+h^{-1}\|{\bf R}^{*}_{h,r}\|_{0,\o_r}\stl\rho_r(h)\|\v^{*}_{h,r}\|_{\ast,\o_r}
\label{stab:5.01newnew}
\e
and
\ee
\|\w^{*}_{h,r}\|_{0,\o_r}+\|p^{\ast}_{h,r}\|_{1,\o_r}\stl\rho_r(h)\|\v^{*}_{h,r}\|_{\#,\o_r}.
\label{stab:5.02newnew}
\e

By (\ref{4.014}), (\ref{decom5.new1}) and (\ref{4.010}), we obtain the decomposition of $\v_h$ on $\o_r\in\Sigma_2$
as \ee \v_{h,r}=\nabla
(p^{*}_{h,r}+\sum\limits_{i\in \Lambda_r^1}\tilde{p}_{h,i})+\r_h(\w^{*}_{h,r}
+\sum\limits_{i\in \Lambda_r^1}\tilde{\w}_{h,i})+{\bf R}^{*}_{h,r}
+\sum\limits_{i\in \Lambda_r^1}\tilde{\bf R}_{h,i}\,. \l{4.015} \e

Let $\tilde{p}^*_{h,r}\in Z_h(\o)$ and $\tilde{\w}^*_{h,r}\in (Z_h(\o))^3$ denote {\it stable} extensions of $p^{\ast}_{h,r}$
and $\w^*_{h,r}$, respectively. Additionally, let $\tilde{\bf R}^{*}_{h,r}\in
V_h(\o)$ denote the standard extension of ${\bf
R}^{*}_{h,r}$ by zero onto $\Omega$. Because $p^{*}_{h,r}$,
$\w^{*}_{h,r}$ and ${\bf R}^{*}_{h,r}$ have zero degrees of freedom on $\hat{\mathcal S}_r$, Lemma \ref{extension-new} implies that the extensions
$\tilde{p}^*_{h,r}$, $\tilde{\w}^*_{h,r}$ and $\tilde{\bf
R}^{*}_{h,r}$ vanish on every $\o_{l}\in \Sigma_1$.
Then, we set \ee \tilde{\v}^{*}_{h,r}=\nabla\tilde{p}^{*}_{h,r}
+\r_h\tilde{\w}^{*}_{h,r}+\tilde{\bf R}^{*}_{h,r}~~ \m{for all $r$
such that} ~~\o_r\in\Sigma_2.\l{4.016} \e

{\bf Step 3}: Obtain the final desired decomposition of $\v_h$.

Assume that the decompositions of $\v_h$ on all polyhedra belonging to $\Sigma_1$,
$\Sigma_2$, $\cdots,\Sigma_{l-1}$ are built. We now consider the index $l\geq 3$, and will build up a
decomposition of $\v_h$ in all subdomains $\o_r\in\Sigma_l$.

As in {\bf Step 2}, we introduce another index set:
$$ \Lambda_r^{l-1}=\{~i~; ~n_1+1\leq i\leq
n_{l-1} \m{~such that
~}\partial\o_i\cap\partial\o_r\not=\emptyset\},\quad\Omega_r\in\Sigma_l~(l\geq3).$$
Define the remainder
\ee
\l{eq:sig2}
\v^{*}_{h,r}=\v_{h,r}-\sum\limits_{i\in\Lambda_r^1}\tilde{\v}_{h,i}
-\sum\limits_{i\in\Lambda_r^{l-1}}\tilde{\v}^{*}_{h,i}\quad\mbox{on}~~\o_r
\e
and
$$ \hat{\mathcal S}_r=\bigcup_{i\in \Lambda_r^1}\bigcup_{i\in\Lambda_r^{l-1}}(\partial\Omega_r\cap\partial\Omega_i)\bigcup(\partial\Omega_r\cap\partial\Omega).  $$

By (\ref{eq:sig2}) and the definitions of $\tilde{\v}_{h,i}$ and $\tilde{\v}^{*}_{h,i}$, we know that $\la_e(\v^{*}_{h,r})=0$ for all
$e\subset\hat{\mathcal S}_r$. Thus, using {\bf Condition B$'$} in Definition \ref{classM} and Theorem 3.1 in Section 3 of \cite{Hu1-2017} again, we find
$p^{*}_{h,r}\in Z_h(\o_r)$, $\w^{*}_{h,r}\in (Z_h(\o_r))^3$ and
${\bf R}^{*}_{h,r}\in V_h(\o_r)$ such that \ee\v^{*}_{h, r}=\nabla
p^{*}_{h,r}+\r_h\w^{*}_{h,r}+{\bf R}^{*}_{h,r}~~
\mbox{on}~~\o_r\,,\l{4.017} \e and they have zero degrees of freedom on $\hat{\mathcal S}_r$.
Moreover, we have
\ee
\|\w^{*}_{h,r}\|_{1,\o_r}+h^{-1}\|{\bf R}^{*}_{h,r}\|_{0,\o_r}\stl\rho_r(h)\|\v^{*}_{h,
r}\|_{\ast,\o_r} \label{stab:5.03newnew}
\e
and
\ee
\|\w^{*}_{h,r}\|_{0,\o_r}+\|p^{\ast}_{h,r}\|_{1,\o_r}\stl\rho_r(h)\|\v^{*}_{h,
r}\|_{\#,\o_r}.\label{stab:5.04newnew}
\e
Using (\ref{eq:sig2}) and (\ref{4.017}), we
have the following decomposition for $\v_h$ on each
$\o_r\in\Sigma_l$ : \beq \v_{h,r}&=&\nabla
(p^{*}_{h,r}+\sum\limits_{i\in\Lambda_r^1}\tilde{p}_{h,i}
+\sum\limits_{i\in\Lambda_r^{l-1}}\tilde{p}^{*}_{h,i})
+~\w^{*}_{h,r}+\sum\limits_{i\in\Lambda_r^1}\tilde{\w}_{h,i}
+\sum\limits_{i\in\Lambda_r^{l-1}}\tilde{\w}^{*}_{h,i}\cr &+&{\bf
R}^{*}_{h,r}+\sum\limits_{i\in\Lambda_r^1}\tilde{\bf R}_{h,i}
+\sum\limits_{i\in\Lambda_r^{l-1}}\tilde{\bf
R}^{*}_{h,i}~~~\mbox{on} ~~\o_r.\l{4.019} \eq

As in Steps 1 and 2, we can extend  $p^{*}_{h,r}$,
$\w^{*}_{h,r}$ and ${\bf R}^{*}_{h,r}$ onto the entire domain $\o$
to obtain $\tilde{p}^{*}_{h,r}$, $\tilde{\w}^{*}_{h,r}$ and $\tilde{\bf
R}^{*}_{h,r}$. Then, by Lemma \ref{extension-new}, the extensions
$\tilde{p}^*_{h,r}$, $\tilde{\w}^*_{h,r}$ and $\tilde{\bf
R}^{*}_{h,r}$ vanish on every $\o_{i}\in \Sigma_j$ for $1\leq j\leq
l-1$. Next, we define \ee \tilde{\v}^{*}_{h,r}=\nabla
\tilde{p}^{*}_{h,r}+\r_h\tilde{\w}^{*}_{h,r}+\tilde{\bf R}^{*}_{h,r}
~~~\m{for all $r$ such that} ~~\o_r\in\Sigma_l. \l{4.021} \e It is clear that
$\la_e(\tilde{\v}^{*}_{h,r})=0$ for all $e\in\hat{\mathcal S}_r$.

Continuing with the above procedure for all $l$ until $l=m$ \footnote{This is a recursive process. Decompositions on the subdomains in
$\Sigma_1$ and $\Sigma_2$ have been built by {\bf Step 1} and {\bf Step 2}, then we establish a decomposition on the subdomains in $\Sigma_3$ by the above procedure (for $l=3$),
and we further build a decomposition on the subdomains in $\Sigma_4$ by the above procedure (for $l=4$), and so on.}, we
build up a decomposition of $\v_h$ over all the
subdomains $\Omega_1$, $\Omega_2$, $\ldots$, $\Omega_{N_0}$ such
that \ee \l{eq:final_decomp}
\v_h=\sum\limits_{r=1}^{n_1}\tilde{\v}_{h,r}
+\sum\limits_{r=n_1+1}^{n_m}\tilde{\v}^{*}_{h,r} = \nabla
p_h+\r_h\w_h+{\bf R}_h \e where $p_h\in Z_h(\o)$, $\w_h\in
(Z_h(\o))^3$ and ${\bf R}_h\in V_h(\o)$ are given
by \ee \l{eq:final} p_h=\sum\limits_{r=1}^{n_1}\tilde{p}_{h,r}
+\sum\limits_{r=n_1+1}^{n_m}\tilde{p}^{*}_{h,r},\quad
\w_h=\sum\limits_{r=1}^{n_1}\tilde{\w}_{h,r}
+\sum\limits_{r=n_1+1}^{n_m}\tilde{\w}^{*}_{h,r}\e and \ee{\bf
R}_h=\sum\limits_{r=1}^{n_1}\tilde{\bf R}_{h,r}
+\sum\limits_{r=n_1+1}^{n_m}\tilde{\bf R}^{*}_{h,r}.\e

\begin{remark} We must emphasize that, for each $\o_r\in \Sigma_l$
($2\leq l\leq m$), the extensions $\tilde{p}^*_{h,r}$,
$\tilde{\w}^*_{h,r}$ and $\tilde{\bf R}^{*}_{h,r}$ vanish on every
$\o_{i}\in \Sigma_j$ for $1\leq j\leq l-1$; otherwise, the
decomposition (\ref{eq:final_decomp}) is not yet valid. This is why
we require $p^*_{h,r}$, $\w^*_{h,r}$ and ${\bf R}^{*}_{h,r}$ to vanish on
$\hat{\mathcal S}_r$, as guaranteed by Theorem 3.1 in Section 3 of  \cite{Hu1-2017}.
Hence, we had to build various regular decompositions in \cite{Hu1-2017}.
\end{remark}

\subsubsection{Stability of the decomposition}
\label{sec:proof}
This part is devoted to the proof of the stability estimates in Theorem~\ref{thm:main1} based on the regular decomposition defined in the previous subsection.

The following results can be obtained directly from the definitions of $\|\v_{h,r}\|^2_{\ast,\o_r}$ and $\|\v_{h,r}\|^2_{\#,\o_r}$.\\
{\bf Proposition 5.1}. Let $\Omega_r\subset\Omega$ and $\v_h\in V_h(\Omega_r)$. For Cases (i)-(iv), we always have
\ee
\al_r\|\v_h\|^2_{\ast,\o_r}\stl\al_r\|\c\,\v_h\|^2_{0,\o_r}+\be_r\|\v_h\|^2_{0,\o_r}
\label{inequality-4.1113}
\e
and
\ee
\be_r\|\v_h\|^2_{\#,\o_r}\stl\al_r\|\c\,\v_h\|^2_{0,\o_r}+\be_r\|\v_h\|^2_{0,\o_r}. \label{inequality-4.1114}
\e
\hfill $\Box$

Let us recall two definitions introduced in Subsection 5.1.2 of \cite{HuShuZou2013}.
\begin{definition}\l{def:4.1}
For a subdomain $\o_r$, another subdomain $\o_{r'}$ is called
a ``parent" of $\o_r$ if $\bar\o_{r'}\cap\bar\o_{r}\not=\emptyset$ and
$\al_{r'} \geq \al_r$.
\end{definition}
\begin{definition}\l{def:4.2} A parent of subdomain $\o_r$ is called a
level-$1$ ancestor of $\o_r$, and a parent of a level-$1$ ancestor
of $\o_r$ is called a level-$2$ ancestor of $\o_r$. In general,
a parent of a level-$j$ ancestor of $\o_r$ is called a
level-$(j+1)$ ancestor of $\o_r$.
\end{definition}

For a polyhedron $\o_r\in\Sigma_l$ ($l\geq 2$),
let $\Lambda_r^{(j)}(a)$ denote the set of all level-$j$ ancestors of $\o_r$.
Note that $\Lambda_r^{(j)}(a)$ may be the empty set for some $j$. For such $\Omega_r$, we use $L_r(a)$ to denote the number of non-empty sets $\Lambda_r^{(j)}(a)$.
Without loss of generality, we assume that the sets $\Lambda_r^{(j)}(a)$ are non-empty for $j=1,\cdots,L_r(a)$.

Let $\rho_r(h)$ be as defined at the beginning of this section, and set $\rho(h)=\max_{1\leq r\leq N_0}\rho_r(h)$. For a subdomain $\Omega_r$, define a norm as
$$ \|\v\|_{H^{\ast}(\c,\Omega_r)}=(\al_r\|\c~\v\|^2_{0,\Omega_r}+\be_r\|\v\|^2_{0,\Omega_r})^{{1\over 2}},\quad\v\in H(\c;\Omega_r). $$
The following auxiliary result can be viewed as an extension of Lemma 6.1 in Section 6 of \cite{HuShuZou2013}.

\begin{lemma}\l{lem:anysub} Assume that {\bf Condition A} in Definition \ref{class} (and Definition \ref{classM}) is satisfied. For a subdomain $\o_r$ from $\Sigma_l$ ($l\geq 2$),
let $\v^{*}_{h,r}$ be as defined in Steps 2 and 3 for the
construction of the decomposition of any $\v_h\in V_h(\o)$ in
Subsection~\ref{sec:decomp}. Then $\v^{*}_{h,r}$ admits the
following estimates
\ee
\al^{{1\over 2}}_r\|\v^{\ast}_{h,r}\|_{\ast,\o_r}\stlm\|\v_h\|_{H^{\ast}(\c,\o_r)}+\sum\limits_{j=1}^{L_r(a)}\rho^{2j}(h)\sum\limits_{i\in
\Lambda_r^{(j)}(a)}\|\v_h\|_{H^{\ast}(\c,\o_i)}\label{inequality-4.1115}
\e
and
\ee
\be^{{1\over 2}}_r\|\v^{\ast}_{h,r}\|_{\#,\o_r}\stlm\|\v_h\|_{H^{\ast}(\c,\o_r)}+\sum\limits_{j=1}^{L_r(a)}\rho^{2j}(h)\sum\limits_{i\in
\Lambda_r^{(j)}(a)}\|\v_h\|_{H^{\ast}(\c,\o_i)}.\label{inequality-4.1116}
\e
In particular, if all of the sets $\hat{\mathcal S}_r$ are connected, the norm $H^{\ast}(\c,\cdot)$ in (\ref{inequality-4.1115}) can be replaced by the corresponding
$\c$ semi-norm.
\end{lemma}

\no{\it Proof}. The proof contains some different details from that of Lemma 6.1 in \cite{HuShuZou2013}. We prove by induction, and start with
the case of $l=2$. It follows from (\ref{decom5.new1}) and (\ref{4.014}) that
\ee
\|\c\,\v^{*}_{h,r}\|_{0,\o_r} \stl\|\c\,\v_{h,
r}\|_{0,\o_r}+\sum\limits_{i\in\Lambda_r^1}(\|\c\,(\r_h\tilde{\w}_{h,i})\|_{0,\o_r}+\|\c~\tilde{{\bf R}}_{h,i}\|_{0,\o_r}).
\l{inequality-4.1111} \e
By {\bf Proposition 4.1} in \cite{Hu1-2017} and Lemma \ref{extension-new}, we deduce that
$$
\|\c(\r_h\tilde{\w}_{h,i})\|_{0,\o_r}=\|\c\tilde{\w}_{h,i}\|_{0,\o_r}\stl\rho_i(h)\|\w_{h,i}\|_{1,\o_{i}}.$$
Moreover, by the inverse estimate and the definition of $\tilde{{\bf R}}_{h,i}$, we have
$$ \|\c~\tilde{{\bf R}}_{h,i}\|_{0,\o_r}\stl h^{-1}\|{\bf R}_{h,i}\|_{0,\o_{i}}. $$
Substituting the above two inequalities into (\ref{inequality-4.1111}) and using (\ref{4.009}) yields
\beqx
\|\c\,\v^{*}_{h,r}\|_{0,\o_r} &\stl&\|\c\,\v_h\|_{0,\o_r}
+\sum\limits_{i\in\Lambda_r^1}\rho_i(h)\|\v_h\|_{\ast,\o_{i}}\cr
&\stl&\|\c\,\v_h\|_{0,\o_r}+\rho(h)\sum\limits_{i\in
\Lambda_r^{(1)}(a)}\|\v_h\|_{\ast,\o_i}. \eqx
Furthermore, by the relation $\al_r\stl\al_i$ for $i\in
\Lambda_r^{(1)}$ and the inequality (\ref{inequality-4.1113}) we obtain
\ee
\al^{{1\over 2}}_r\|\c\,\v^{*}_{h,r}\|_{0,\o_r}\stl\al^{{1\over 2}}_r\|\c\,\v_h\|_{0,\o_r}+\rho(h)\sum\limits_{i\in
\Lambda_r^{(1)}(a)}\|\v_h\|_{H^{\ast}(\c,\o_i)}.\label{inequality-4.1112}
\e
Similarly, using (\ref{stab:5.05newnew}) and (\ref{inequality-4.1114}), together with {\bf Condition A} in Definition \ref{class} (and Definition \ref{classM}), we can show that
$$ \be^{{1\over 2}}_r\|\v^{*}_{h,r}\|_{0,\o_r}\stlm\be^{{1\over 2}}_r\|\v_h\|_{0,\o_r}+\rho(h)\sum\limits_{i\in
\Lambda_r^{(1)}(a)}\|\v_h\|_{H^{\ast}(\c,\o_i)}.$$
Combining the above inequality with (\ref{inequality-4.1112}), we know that
(\ref{inequality-4.1115}) and (\ref{inequality-4.1116}) are valid for all the subdomains $\o_r$ in
$\Sigma_2$.

Now, assume that (\ref{inequality-4.1115}) is true for all subdomains
$\o_r\in\Sigma_l$ with $l\leq n$. We need to verify
(\ref{inequality-4.1115}) for all subdomains $\o_r\in\Sigma_{n+1}$.
It follows from (\ref{eq:sig2}) that
\beq \|\c\,\v^{*}_{h,r}\|_{0,\o_r} &\stl&
\|\c\,\v_{h,r}\|_{0,\o_r}
+\sum\limits_{i\in\Lambda_r^1}\|\c\,\tilde{\v}_{h,i}\|_{0,\o_r}\cr
&+&\sum\limits_{i\in\Lambda_r^{n}}\|\c\,\tilde{\v}^{*}_{h,i}\|_{0,\o_r}.
\l{4.new2}
\eq
Similarly as (\ref{inequality-4.1112}) was derived, one can check that (by (\ref{4.016}), Lemma \ref{extension-new} and (\ref{stab:5.01newnew}))
\beqx
\|\c\,\tilde{\v}_{h,i}\|_{0,\o_r}&\stl&\rho_i(h)(\|\w_{h,i}\|_{1,\o_i}+h^{-1}\|{\bf R}_{h,i}\|_{0,\o_i})\stl\rho^2_i(h)
\|\v_{h,i}\|_{\ast,\o_i}~~~(i\in\Lambda_r^1),\\
\|\c\,\tilde{\v}^{*}_{h,i}\|_{0,\o_r}&\stl&
\rho_i(h)(\|\w^{\ast}_{h,i}\|_{1,\o_i}+h^{-1}\|{\bf R}^{\ast}_{h,i}\|_{0,\o_i})\stl\rho^2_i(h)
\|\v^{\ast}_{h,i}\|_{\ast,\o_i}~~~(i\in\Lambda_r^n).
\eqx
Combining these estimates with (\ref{4.new2}) gives
\beq\al_r^{{1\over 2}}\|\c\,\v^{*}_{h,r}\|_{0,\o_r}&\stl&
\al_r^{{1\over 2}}\|\c\,\v_{h,r}\|_{0,\o_r}+\rho^2(h)
\sum\limits_{i\in\Lambda_r^1}\al_i^{{1\over 2}}\|\v_{h,i}\|_{\ast,\o_i}\cr
&&+~\rho^2(h)\sum\limits_{i\in\Lambda_r^{n}}\al_i^{{1\over 2}}\|\v^{*}_{h,i}\|_{\ast,\o_i}.
\l{4.new3}\eq
Here we have used the fact that $\al_r\stl\al_i$ for $i\in\Lambda_r^1\cup\Lambda_r^{n}$.
Noting that for $i\in\Lambda_r^{n}$, we have
$\o_i\in\Sigma_l$ for some $l\leq n$. Thus, by the inductive
assumption,
\beq \l{eq:induct}
\sum\limits_{i\in\Lambda_r^{n}}\al^{{1\over 2}}_i\|\v^{*}_{h,i}\|_{\ast,\o_i}&\stl&
\sum\limits_{i\in\Lambda_r^{n}}\|\v_h\|_{H^{\ast}(\c,\o_i)}\cr
&+&\sum\limits_{i\in\Lambda_r^{n}}
\sum\limits_{j=1}^{L_i(a)}\rho^{2j}(h)\sum\limits_{k\in
\Lambda_i^{(j)}(a)}\|\v_h\|_{H^{\ast}(\c,\o_k)}. \eq But, for all
subdomains $\o_r\in\Sigma_{n+1}$ and $i\in\Lambda_r^{n}$, we know
$L_i(a)\le L_r(a)$ and $\Lambda_i^{(j)}(a)=\emptyset$ for $j>L_i(a)$
by definition, so we have the relation
\ee
\sum\limits_{i\in\Lambda_r^{n}}\sum\limits_{j=1}^{L_i(a)}\sum\limits_{k\in
\Lambda_i^{(j)}(a)}
=\sum\limits_{j=1}^{L_r(a)}\sum\limits_{i\in\Lambda_r^{n}}\sum\limits_{k\in
\Lambda_i^{(j)}(a)}
=\sum\limits_{j=1}^{L_r(a)} \sum\limits_{k\in \Lambda_r^{(j+1)}(a)}.
\l{relation}\e
Combining this with the fact that
$\Lambda_r^{(j+1)}(a)=\emptyset$ for $j\geq L_r(a)$, we get
$$
\sum\limits_{i\in\Lambda_r^{n}}\sum\limits_{j=1}^{L_i(a)}\sum\limits_{k\in
\Lambda_i^{(j)}(a)}=\sum\limits_{j=1}^{L_r(a)-1} \sum\limits_{k\in
\Lambda_r^{(j+1)}(a)}.
$$
From this identity and (\ref{eq:induct}) it follows that
\beq \l{eq:induct1}
\sum\limits_{i\in\Lambda_r^{n}}\al^{{1\over 2}}_i\|\v^{*}_{h,i}\|_{\ast,\o_i}&\stl&
\sum\limits_{i\in\Lambda_r^{n}}\|\v_h\|_{H^{\ast}(\c,\o_i)}\cr
&+&\sum\limits_{j=2}^{L_r(a)}\rho^{2j}(h)\sum\limits_{k\in
\Lambda_r^{(j)}(a)}\|\v_h\|_{H^{\ast}(\c,\o_k)}.
\eq
Substituting this into (\ref{4.new3}), and using (\ref{inequality-4.1113}) and the identity
$$ \sum\limits_{i\in\Lambda_r^1}+\sum\limits_{i\in \Lambda_r^n}
=\sum\limits_{i\in\Lambda_r^{(1)}(a)}\quad \m{for} ~~\o_r\in\Sigma_{n+1},$$
we can
immediately derive that \ee
\al^{{1\over 2}}_r\|\c\,\v^{*}_{h,r}\|_{0,\o_r}\stl \al^{{1\over 2}}_r\|\c\,\v_h\|_{0,\o_r}
+\sum\limits_{j=1}^{L_r(a)}\rho^{2j}(h)\sum\limits_{k\in
\Lambda_r^{(j)}(a)}\|\v_h\|_{H^{\ast}(\c,\o_k)}.\l{4.new4} \e
Similarly, by using (\ref{stab:5.02newnew}) and (\ref{inequality-4.1114}), and noticing {\bf Condition A} in Definition \ref{class} (and Definition \ref{classM}), we can show that
$$ \be^{{1\over 2}}_r\|\v^{*}_{h,r}\|_{0,\o_r}\stlm \be^{{1\over 2}}_r\|\v_h\|_{0,\o_r}
+\sum\limits_{j=1}^{L_r(a)}\rho^{2j}(h)\sum\limits_{k\in
\Lambda_r^{(j)}(a)}\|\v_h\|_{H^{\ast}(\c,\o_k)}. $$
Combining the above inequality with (\ref{4.new4}), we prove (\ref{inequality-4.1115}) and (\ref{inequality-4.1116}) for all subdomains
$\o_r\in\Sigma_{n+1}$, thus we complete the proof of the estimates by the mathematical induction.

The final conclusion can be seen directly from the above proof (refer to the proof of Lemma 6.1 in \cite{HuShuZou2013}).
\hfill $\Box$

\medskip

\no{Proof of Theorem~\ref{thm:main1}}. We are now ready to prove Theorem~\ref{thm:main1}. We start with the estimate of $\|\w_h\|^2_{H^1_{\omega}(\o_r)}$ for
each subdomain $\o_r$ in $\Sigma_1$, i.e., $1\le r\leq n_1$. For this case, we have $\w_h|_{\o_r}=\w_{h,r}$ with $\w_{h,r}$ being
defined by (\ref{4.001}) (note that any two of the subdomains $\o_1,\cdots,\o_{n_1}$ do not intersect). Then, it follows from (\ref{4.009}) and {\bf Proposition 5.1} that
\ee\|\w_h\|^2_{H^1_{\omega}(\o_r)}\stl
\al_r\|\v_{h,r}\|^2_{\ast,\o_r}\stl\al_r\|\c\,\v_h\|^2_{0,\o_r}+\be_r\|\v_h\|^2_{0,\o_r}.
\l{4.023}\e

Next, we consider all the subdomains $\o_r$ in $\Sigma_2$. From the construction in Step 2, we have
$$
\w_h|_{\o_r}=\tilde{\w}^{*}_{h,r}+(\sum\limits_{i\in\Lambda_r^1}\tilde{\w}_{h,i})|_{\o_r}.$$
It follows by the triangle inequality that
\ee
\|\w_h\|_{1,\o_r}\stl\|\tilde{\w}^{*}_{h,r}\|_{1,\o_r}
+\sum\limits_{i\in\Lambda_r^1}\|\tilde{\w}_{h,i}\|_{1,\o_r}.\l{4.024}
\e
Because $\tilde{\w}^{*}_{h,r}=\w^{*}_{h,r}$ on $\o_r$, we obtain the following inequality from (\ref{stab:5.01newnew})
$$
\|\tilde{\w}^{*}_{h,r}\|_{1,\o_r}\stl\rho_r(h)\|\v^{*}_{h,r}\|_{\ast,\o_r}.
$$
Substituting this and (\ref{4.009}) (noting the stability of the extension) into (\ref{4.024}) leads to
\ee
\|\w_h\|_{1,\o_r}\stl\rho_r(h)\|\v^{*}_{h,r}\|_{\ast,\o_r}
+\sum\limits_{i\in\Lambda_r^1}\|\v_{h,r}\|_{\ast,\o_i}.\l{4.027}
\e
Then, by inserting the coefficient $\omega$ and using (\ref{inequality-4.1115}) for $l=2$, we readily have that, for
all subdomains $\o_r\in\Sigma_2$,
\beq
\al^{{1\over 2}}_r\|\w_h\|_{1,\o_r}&\stlm&\rho_r(h)\al^{{1\over 2}}_r\|\v^{\ast}_{h,r}\|_{\ast,\o_r}
+\sum\limits_{i\in\Lambda_r^1}\al_r^{{1\over 2}}\|\v_{h,i}\|_{\ast,\o_i}\cr
&\stlm&\rho(h)\|\v_h\|_{H^{\ast}(\c,\o_r)}
+\rho^3(h)\sum\limits_{i\in\Lambda_r^1}\|\v_h\|_{H^{\ast}(\c,\o_i)}. \l{4.025} \eq
Here, we have used the fact that $\al_r\leq \al_i$ for $i\in\Lambda_r^1$, which comes from the definitions of $\Sigma_1$ and $\Sigma_2$.

Finally, we consider all the subdomains $\o_r$ from the general class
$\Sigma_l$ with $l\geq 3$. {By the definition of $\w_h$ (see
(\ref{eq:final})), we have that, for $\w_h$ in $\o_r$,
\ee
\w_h=\w^{\ast}_{h,r}
+\sum\limits_{i\in\Lambda_r^1}\tilde{\w}_{h,i}
+\sum\limits_{i\in\Lambda_r^{l-1}}\tilde{\w}^{*}_{h,i}.\label{equality:5.1}
\e
Then,
\ee
\|\w_h\|_{1,\o_r}\stl\|\w^{\ast}_{h,r}\|_{1,\o_r}
+\sum\limits_{i\in\Lambda_r^1}\|\tilde{\w}_{h,i}\|_{1,\o_r}
+\sum\limits_{i\in\Lambda_r^{l-1}}\|\tilde{\w}^{*}_{h,i}\|_{1,\o_r}.\label{equal:5.01newnew}
\e
By Lemma \ref{extension-new}, we have
$$ \|\tilde{\w}_{h,i}\|_{1,\o_r}\stl\rho_i(h)\|\w_{h,i}\|_{1,\o_i}~~~(i\in\Lambda_r^1) $$
and
$$ \|\tilde{\w}^{\ast}_{h,i}\|_{1,\o_r}\stl\rho_i(h)\|\w^{\ast}_{h,i}\|_{1,\o_i}~~~(i\in\Lambda_r^{l-1}). $$
Together with (\ref{4.009}) and (\ref{stab:5.01newnew}), this leads to
$$ \|\tilde{\w}_{h,i}\|_{1,\o_r}\stl\rho_i(h)\|\v_h\|_{\ast,\o_i}~~~(i\in\Lambda_r^1) $$
and
$$ \|\tilde{\w}^{\ast}_{h,i}\|_{1,\o_r}\stl\rho^2_i(h)\|\v^{\ast}_{h,i}\|_{\ast,\o_i}~~~(i\in\Lambda_r^{l-1}). $$
Substituting (\ref{stab:5.03newnew}) and the above two inequalities into (\ref{equal:5.01newnew}) yields
\beq
\|\w_h\|_{1,\o_r}\stl\rho_r(h)\|\v^{\ast}_{h,r}\|_{\ast,\o_r}
+\sum\limits_{i\in\Lambda_r^1}\rho_i(h)\|\v_h\|_{\ast,\o_i}
+\sum\limits_{i\in\Lambda_r^{l-1}}\rho^2_i(h)\|\v^{*}_{h,i}\|_{\ast,\o_i}.\label{equal:5.02newnew}
\eq
Inserting the coefficient $\omega_r$ into the above inequality and using (\ref{inequality-4.1115}) and the relation $\al_r\leq \al_i$ (for $i\in \Lambda_r^{(j)}(a)$) gives
\begin{eqnarray*}
\al^{\12}_r\|\w_h\|_{1,\o_r}
\stlm\rho_r(h)\|\v_h\|_{H^{\ast}(\c,\o_r)}+\rho_r(h)\sum\limits_{j=1}^{L_r(a)}\rho^{2j}(h)\sum\limits_{i\in
\Lambda_r^{(j)}(a)}\|\v_h\|_{H^{\ast}(\c,\o_i)}.
\end{eqnarray*}
It is clear that the set $\Lambda_r^{(j)}(a)$ contains only a few indices $i$ and that $L_r(a)$ is a finite number.
Combining the above estimate with those in (\ref{4.023}) and
(\ref{4.025}), we obtain \beq
\sum_{r=1}^{N_0}\al_r\|\w_h\|^2_{1,\o_r}&\stlm&\sum_{r=1}^{N_0}(\rho^2_r(h)\|\v_h\|^2_{H^{\ast}(\c,\o_r)})\cr
&&+~\sum\limits_{r={n_1+1}}^{N_0}\rho_r(h)\sum\limits_{j=1}^{L_r(a)}
\rho^{4j}(h)\sum\limits_{i\in \Lambda_r^{(j)}(a)}\|\v_h\|^2_{H^{\ast}(\c,\o_i)}\cr
&&\stlm\rho^{2m}(h)\sum\limits_{r=1}^{N_0}
\|\v_h\|^2_{H^{\ast}(\c,\o_r)},\l{4.new6} \eq
where $m=\max\limits_{1\leq r\leq N_0}(2L_r(a)+1)$. It follows from (\ref{4.new6}) that
$$ (\sum_{r=1}^{N_0}\al_r\|\w_h\|^2_{1,\o_r})^{\12}\stlm\rho^{m}(h)\|\v_h\|_{H^{\ast}(\c,\o)}. $$
In an analogous way, but using (\ref{inequality-4.1116}) and {\bf Condition A} in Definition \ref{class}, we can verify that
$$ \|\be^{{1\over 2}}\w_h\|_{0,\Omega},~~\|p_h\|_{H^1_{\beta}(\Omega)}\stlm\rho^{m}(h)\|\v_h\|_{H^{\ast}(\c,\o)}.$$
The term
$\|{\bf R}_h\|^2_{L^2_{\omega}(\o_r)}$ can be estimated more easily. This completes the proof of Theorem~\ref{thm:main1}.
\hfill $\Box$


\subsection{Analysis for the case with thorny vertices}

In this subsection, we present a proof of Theorem \ref{thm:main}. Here we assume that the coefficients $\{\al_k\}$ and $\{\be_k\}$ belong to the problem class $\mathbb{P}(\mu)$.
To make the analysis easier to understand, we first describe the basic idea of the proof.

For a {\it thorny vertex} $\vv$, let $\Im^{\ast}_{\vv}$ be the set defined in Definition \ref{definition3}, and let $\Im^{\ast}$ and $\Im^c_{\ast}$ be defined in Subsection 3.2, i.e.,
$$\Im^{\ast}=\cup_{\vv\in {\mathcal V}_{\ast}}\Im^{\ast}_{\vv},\quad\Im^c_{\ast}=\{\Omega_k\}_{k=1}^{N_0}\backslash\Im^{\ast}.$$


If there is no {\it thorny vertex} (i.e., $\Im^{\ast}=\emptyset$, which implies that $\Im^c_{\ast}=\{\Omega_k\}_{k=1}^{N_0}$), then the distribution of the coefficients satisfy the
{\it generalized quasi-monotonicity assumption} by {\bf Proposition 3.4}.
Then, the desired regular decomposition can be directly obtained by Theorem \ref{thm:main1}. This inspires us to separately construct regular decompositions on the polyhedra
belonging to $\Im^{\ast}$ and $\Im^c_{\ast}$. Unfortunately, building the desired regular decomposition on the polyhedra in $\Im^{\ast}$ is somewhat technical, because the union of
all these polyhedra is a non-Lipschitz domain. For convenience, each polyhedron in $\Im^{\ast}$ is called a {\it thorny polyhedron}. Theorems 3.3 and 3.4 in Section 3 of \cite{Hu1-2017}
provide the basic tool for building the desired regular decompositions on the {\it thorny polyhedra}. However, as a {\it thorny polyhedron} may have
several {\it thorny vertices}, such a decomposition cannot be directly obtained. Hence, we have to construct auxiliary subdomains, each of which contains only a {\it thorny vertex}.
\subsubsection{Auxiliary subdomains associated with thorny vertices}
For $\vv_i\in {\mathcal V}_{\ast}$, let $\{\Omega_{i_r}\}_{r=1}^{m_i}$ denote all the polyhedra in $\Im^{\ast}_{\vv_i}$. For each polyhedron $\Omega_{i_r}$, we cut its angle containing $\vv_i$ as
in Example 4.2 in Subsection 4.1 of \cite{Hu1-2017}, and use $D_{i_r}$ to denote the cut domain (see the left one in Figure 6 and compare it with Figure 2).

\begin{center}
\includegraphics[width=6cm,height=4.3cm]{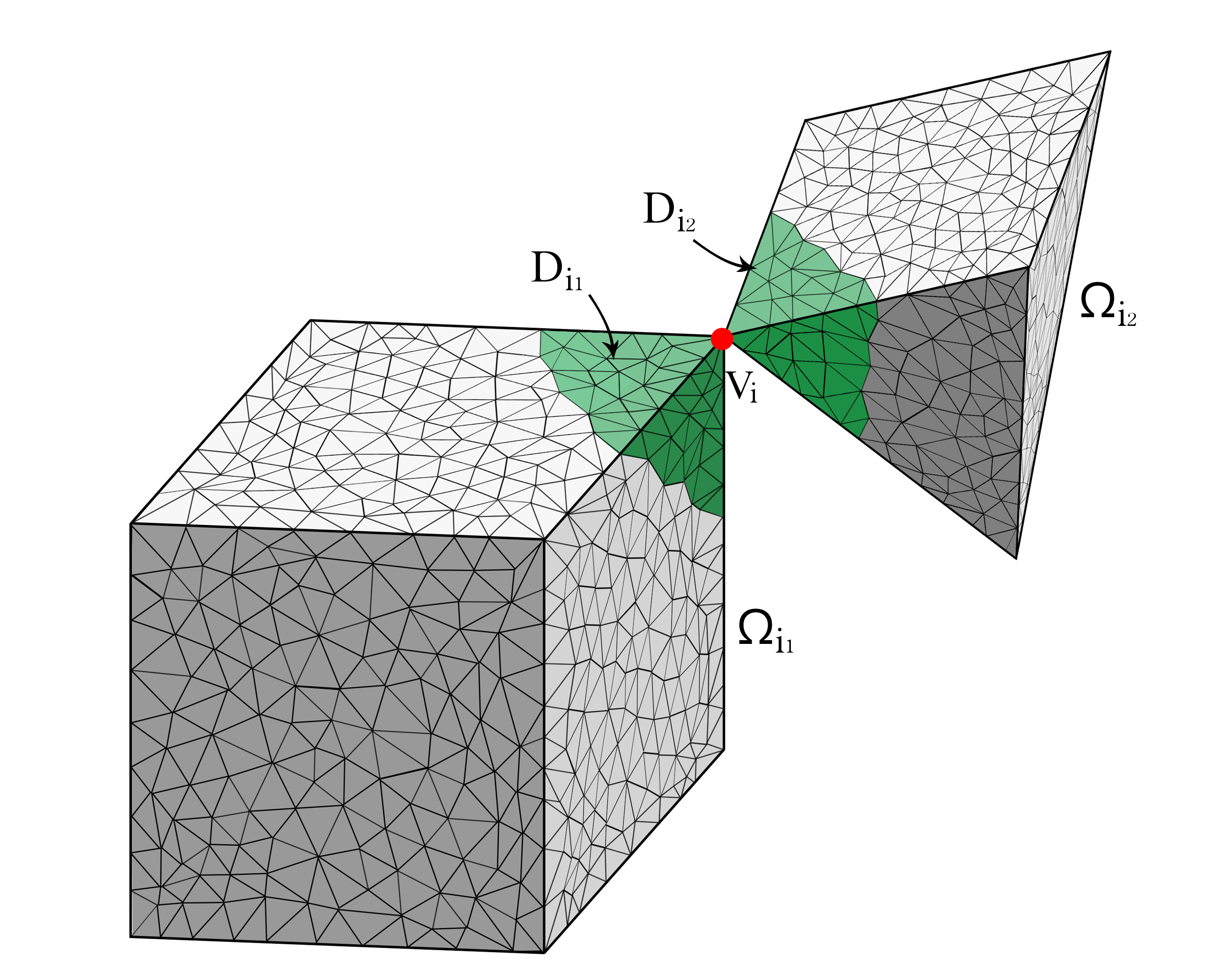}\quad
\includegraphics[width=5.5cm,height=4.1cm]{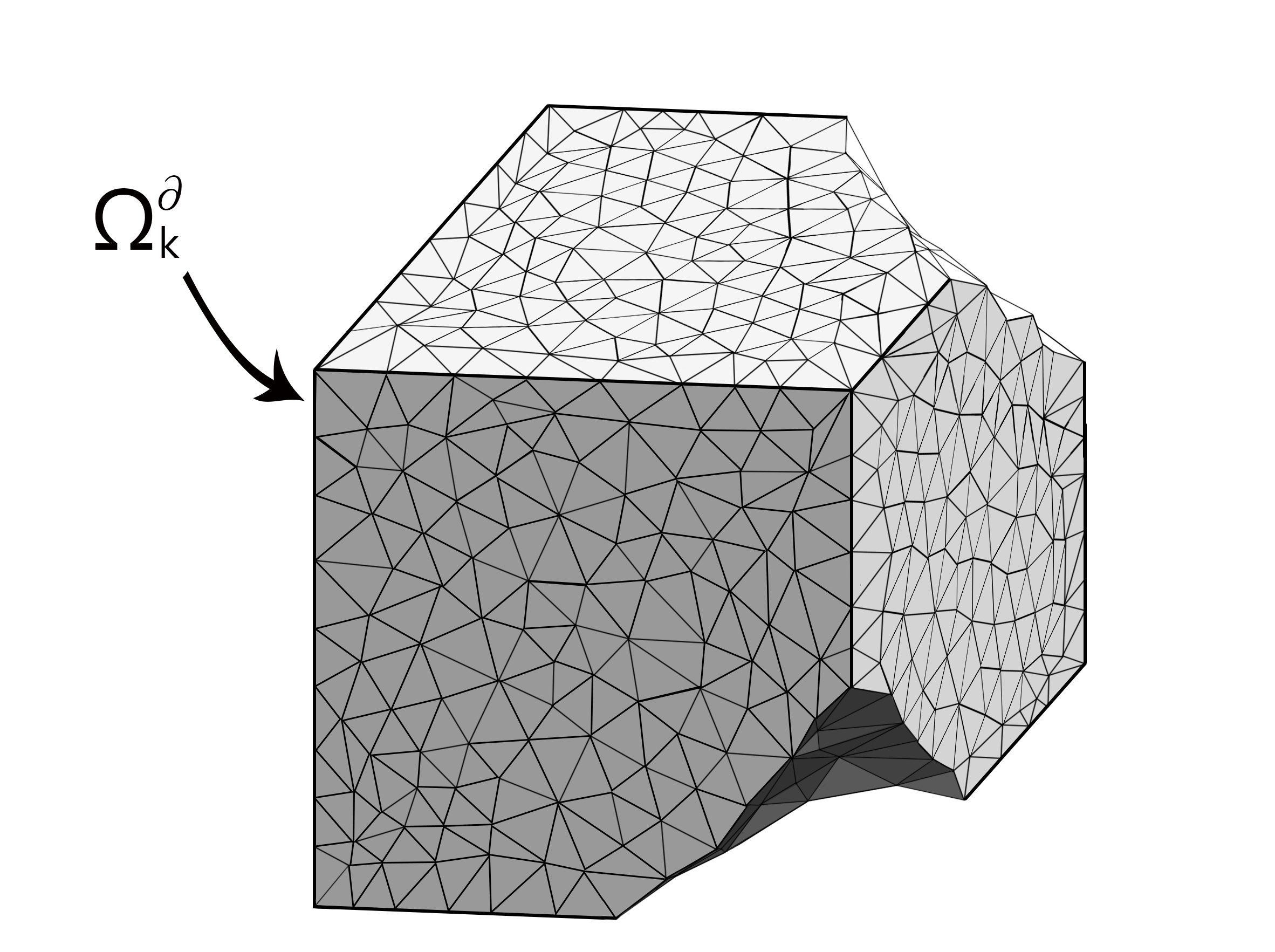}

\centerline{Figure 6. left: The green part denote {\it cut domains};~~ right: {\it residual domain} from $\Omega_k$}
\end{center}
\vskip 2mm

The cut domains satisfy the following conditions:
(i) each cut domain $D_{i_r}$ is the union of some tetrahedron elements contained in $\Omega_{i_r}$; (ii) each $D_{i_r}$ is a Lipschitz polyhedron with a Lipschitz constant that is independent of $h$
(it has a size of $O(1)$ and can be regarded as a small perturbation of a usual polyhedron); and (iii) for two different {\it thorny vertices}
$\vv_i$ and $\vv_j$, the cut domains $D_{i_r}$ and $D_{j_l}$ do not intersect each other, i.e., $\bar{D}_{i_r}\cap\bar{D}_{j_l}=\emptyset$. Set
$$
D^{\ast}_{\vv_i}=\bigcup_{r=1}^{m_i}D_{i_r},\quad \Omega^{\ast}_{\vv_i}=\bigcup_{r=1}^{m_i}\Omega_{i_r}.$$
It is clear that $D^{\ast}_{\vv_i}$ is a real subset of $\Omega^{\ast}_{\vv_i}$. Moreover, both $\Omega^{\ast}_{\vv_i}$ and $D^{\ast}_{\vv_i}$ are non-Lipschitz domains like the one
considered in Theorem 3.4 in Section 3 of \cite{Hu1-2017}.

Let $\Lambda^{\ast}$ denote the set of indices $k$ such that $\Omega_k$ has at least one {\it thorny vertex}.
Define the {\it residual domains} of the cut domains as (see the right one in Figure 6)
$$ \Omega_k^{\partial}=\Omega_k\backslash \bigcup_{D_{i_r}\subset \Omega_k}D_{i_r}~~~(k\in\Lambda^{\ast}),\quad D^{\partial}=\bigcup_{k\in \Lambda^{\ast}}\Omega_k^{\partial}. $$
Then, every $\Omega_k^{\partial}$ is a Lipschitz polyhedron with a Lipschitz constant that is independent of $h$. Define
$$ D^{\ast}=\bigcup_{\vv_i\in{\mathcal V}_{\ast}}D^{\ast}_{\vv_i},\quad\Omega^c=\bigcup_{\vv_i\in{\mathcal V}_{\ast}}\bigcup_{\Omega_k\in\Im^c_{\vv_i}}\Omega_k,\quad\Omega^c_{\ast}=\bigcup_{\Omega_k\in\Im^c_{\ast}}\Omega_k,\quad\mbox{and}\quad\Omega^{\partial}=D^{\partial}\cup\Omega^c_{\ast}.$$
Then, the domain $\Omega$ can be decomposed into
$$ \Omega=D^{\ast}\cup D^{\partial}\cup\Omega^c_{\ast}=D^{\ast}\cup\Omega^{\partial}. $$
The main idea is to construct the desired regular decompositions on $D^{\ast}$ and $\Omega^{\partial}$ in turn, where the decomposition on $D^{\ast}$ is defined by
regular decompositions on $\Omega^{\ast}_{\vv_i}$ for all {\it thorny vertices} $\vv_i$.

An extension result associated with each $D^{\ast}_{\vv_i}$ will play a key role in the proof of Theorem \ref{thm:main}.
\begin{lemma}\label{extension1} For a subdomain $D^{\ast}_{\vv_i}$, there exists an extension operator $E_{i,h}$ mapping $Z_h(D^{\ast}_{\vv_i})$ into $Z_h(\Omega)$ such that,
for any function $\phi_h\in Z_h(D^{\ast}_{\vv_i})$, the function $E_{i,h}\phi_h$ satisfies the following conditions:(i) $E_{i,h}\phi_h=\phi_{h}$ on $\bar{D}^{\ast}_{\vv_i}$;
(ii) for any subdomain $G\subset\Omega\backslash D^{\ast}_{\vv_i}$, the function $E_{i,h}\phi_h$ vanishes at all nodes on $\partial G\backslash\bar{D}^{\ast}_{\vv_i}$;
(iii) $E_{i,h}\phi_h$ is discrete harmonic in every subdomain $G\subset\Omega\backslash D^{\ast}_{\vv_i}$ and has the stability estimate
\ee
\|E_{i,h}\phi_h\|_{1,\Omega}\leq C\log (1/h)\|\phi_h\|_{1,D^{\ast}_{\vv_i}}.
\label{stability.new}
\e
\end{lemma}

The proof of this lemma is essentially contained in the proof of Theorem 4.1 in Subsection 4.2 of \cite{Hu1-2017}: the extension $E_{i,h}\phi_h$ can be constructed and analyzed as the extension
$\tilde{p}_{h,1}$ constructed there.

\subsubsection{Proof of Theorem \ref{thm:main}} We can now build a regular decomposition of $\v_h\in V_h^{\ast}(\Omega)$ and establish the corresponding stability estimates.
The proof is divided into three steps.\\
\medskip
{\bf Step 1}. Build regular decompositions on the auxiliary subdomain $D^{\ast}$.

Because $\v_h\in V_h^{\ast}(\Omega)$, the definition of $V_h^{\ast}(\Omega)$ (see Subsection 3.2) implies that the function $\v_h$ satisfies the key assumptions in Theorem 3.4
in Section 3 of \cite{Hu1-2017} for each subdomain $\Omega^{\ast}_{\vv_i}$. Therefore, we can build a regular decomposition of $\v_h$ on every $\Omega^{\ast}_{\vv_i}$, i.e.,
\ee
\v_h=\nabla p^{\ast}_{\vv_i,h}+\r_h\w^{\ast}_{\vv_i,h}+{\bf R}^{\ast}_{\vv_i,h},\quad \mbox{on}~\Omega^{\ast}_{\vv_i}, \label{decom-5.1new}
\e
where $p^{\ast}_{\vv_i,h}\in Z_h(\Omega^{\ast}_{\vv_i})$, $\w^{\ast}_{\vv_i,h}\in (Z_h(\Omega^{\ast}_{\vv_i}))^3$ and ${\bf R}^{\ast}_{\vv_i,h}\in V_h(\Omega^{\ast}_{\vv_i})$. Moreover,
the following estimates hold:
  \ee
  \|p^{\ast}_{\vv_i,h}\|_{1, \Omega_{i_r}}\stl  \log(1/h)\|\v_h\|_{\c, \Omega_{i_r}},\quad\forall \Omega_{i_r}\subset \Omega^{\ast}_{\vv_i}\label{stab-5.1new} \e
and
  \ee \|\w^{\ast}_{\vv_i,h}\|_{1, \Omega_{i_r}}+h^{-1}\|{\bf R}^{\ast}_{\vv_i,h}\|_{0,\Omega_{i_r}}\stl\log(1/h)\|\c~\v_h\|_{0,\Omega_{i_r}},\quad\forall\Omega_{i_r}\subset \Omega^{\ast}_{\vv_i}.\label{stab-5.2new}\e

As $D^{\ast}_{\vv_i}\subset\Omega^{\ast}_{\vv_i}$ and $\bar{D}^{\ast}_{\vv_i}\cap \bar{D}^{\ast}_{\vv_j}=\emptyset$ for $i\not=j$, we can naturally define $p^{\ast}_h\in Z_h(D^{\ast})$,
$\w^{\ast}_h\in (Z_h(D^{\ast}))^3$ and ${\bf R}^{\ast}_h\in V_h(D^{\ast})$ as follows:
$$ p^{\ast}_h=p^{\ast}_{\vv_i,h},~~\w^{\ast}_h=\w^{\ast}_{\vv_i,h},~~{\bf R}^{\ast}_{h}={\bf R}^{\ast}_{\vv_i,h}~~~\mbox{on}~~D^{\ast}_{\vv_i}~~(\forall\vv_i\in{\mathcal V}_{\ast}). $$
Then, we have the decomposition
\ee \v_h=\nabla p^{\ast}_{h}+\r_h\w^{\ast}_{h}+{\bf R}^{\ast}_{h},\quad \mbox{on}~D^{\ast}.\label{decom-5.2new}
\e
{\bf Step 2}. Build a regular decomposition on the global domain $\Omega$.

Let $E_{i,h}$ be the extension operator given in Lemma \ref{extension1}. For the function $p_h^{\ast}$, we define its extension $\tilde{p}_{h}^{\ast}\in Z_h(\Omega)$ as
$\tilde{p}_{h}^{\ast}=\sum_{\vv_i\in{\mathcal V}_{\ast}}E_{i,h}p_{\vv_i,h}^{\ast}$ (here we have not used the extension given by Lemma \ref{extension-new}
because $D^{\ast}_{\vv_i}$ is a non-Lipschitz domain).
Similarly, we can define an extension $\tilde{\w}^{\ast}_{h}$ of $\w^{\ast}_{h}$ such that $\tilde{\w}_{h}^{\ast}\in (Z_h(\Omega))^3$.
Let $\tilde{\bf R}^{\ast}_{h}\in V_h(\Omega)$ be the natural zero extension of ${\bf R}^{\ast}_{h}$. Define the remainder
\ee
\v^{\partial}_h=\v_h-\big(\nabla \tilde{p}^{\ast}_{h}+\r_h\tilde{\w}_{h}^{\ast}+\tilde{\bf R}_{h}^{\ast}\big)\quad\mbox{on}~~\Omega.\label{decom-5.6new}
\e
It follows from (\ref{decom-5.2new}) that $\v_h^{\partial}$ vanishes on $\bar{D}^{\ast}$. By the definition of the {\it cut domains}, there is no {\it thorny vertex} on
the residual domain $D^{\partial}$. Moreover, by {\bf Proposition 3.2}, the domain $\Omega_{\ast}^c$ does not contain
any {\it thorny vertex}, and so there is no {\it thorny vertex} on the domain $\Omega^{\partial}$.
Therefore, from the assumptions given in Theorem \ref{thm:main}, we know that {\bf Condition A} and {\bf Condition B} in the definition of the class $\mathbb{P}(\gamma)$ are satisfied for all
subdomains in $\Omega^{\partial}$. Note that the domain $\Omega^{\partial}$ is a union of Lipschitz polyhedra with Lipschitz constants that are independent of $h$,
so we can apply Theorem \ref{thm:main1} to build a regular decomposition of $\v^{\partial}_h$ on $\Omega^{\partial}$\footnote{The main results established in \cite{Hu1-2017} still hold when
the considered polyhedra are replaced by general Lipschitz polyhedra with Lipschitz constants independent of $h$, where the ``edge lemma" and ``face lemma" needed to be used can be
proved as in Subsection 5.2 and Subsection 5.3 of \cite{Bringmans2020} by using the discrete norms}:
\ee
\v^{\partial}_h=\nabla p^{\partial}_h+\r_h\w^{\partial}_h+{\bf R}^{\partial}_h~~~\mbox{on}~~\Omega^{\partial},\label{decom-5.7new}
\e
where $p^{\partial}_h$, $\w^{\partial}_h$ and ${\bf R}^{\partial}_h$ have zero degrees of freedom on $\bar{D}^{\ast}\cap\bar{\Omega}^{\partial}$. Then, we can naturally define their zero extension
functions $\tilde{p}^{\partial}_h\in Z_h(\Omega)$, $\tilde{\w}^{\partial}_h\in (Z_h(\Omega))^3$ and $\tilde{\bf R}^{\partial}_h\in V_h(\Omega)$, and we have the following decomposition by (\ref{decom-5.7new})
(because $\v_h^{\partial}$ vanishes on $\bar{D}^{\ast}$):
\ee
\v^{\partial}_h=\nabla \tilde{p}^{\partial}_h+\r_h\tilde{\w}^{\partial}_h+\tilde{\bf R}^{\partial}_h~~~\mbox{on}~~\Omega. \label{decom-5.8new}
\e
Moreover, using Theorem \ref{thm:main1} on $\Omega^{\partial}$ and noting that $\tilde{p}^{\partial}_h$, $\tilde{\w}^{\partial}_h$ and $\tilde{\bf R}^{\partial}_h$ vanish on $\bar{D}^{\ast}$,
we obtain the following estimates
 \ee
\|\tilde{p}^{\partial}_h\|_{H^1_{\beta}(\Omega)}\stlm
  \rho^m(h)\|\v^{\partial}_h\|_{H^{\ast}(\c,\Omega^{\partial})}\label{stab-5.3new} \e
and
 \ee \|\tilde{\w}^{\partial}_h\|_{H^1_{\ast}(\o)}+h^{-1}\|\tilde{\bf R}^{\partial}_h\|_{L^2_{\omega}(\Omega)}\stlm
  \rho^m(h)\|\v^{\partial}_h\|_{H^{\ast}(\c,\Omega^{\partial})}.\label{stab-5.4new}\e
Here the exponent $m$ is generally smaller than the exponent $m_0$ in Theorem \ref{thm:main1} because the number of subdomains $\Omega_k\subset\Omega^{\partial}$
is less than $N_0$.

Define
$$ p_h=\tilde{p}^{\ast}_{h}+\tilde{p}^{\partial}_h,\quad \w_h=\tilde{\w}_{h}^{\ast}+\tilde{\w}^{\partial}_h,\quad{\bf R}_h=\tilde{\bf R}_{h}^{\ast}+\tilde{\bf R}^{\partial}_h. $$
By (\ref{decom-5.6new}) and (\ref{decom-5.7new}), we obtain the final decomposition
\ee
\v_h=\nabla p_h+\r_h\w_h+{\bf R}_h\quad\mbox{on}~~\Omega. \label{final-decom}
\e
{\bf Step 3}. Verify the weighted stability of the regular decomposition.

Let $G$ denote a generic subdomain (Lipschitz polyhedron) contained in $\Omega^{\partial}$, and let $\al_G$ and $\be_G$ denote the restrictions of the coefficients on $G$.
We use $\Lambda^{\ast}_G$ to denote the set of indices $i_{r}$ for which the corresponding cut domain $D_{i_r}\subset D^{\ast}_{\vv_i}$ has a common face or common
edge with $G$. By definition, the extension function $\tilde{p}^{\ast}_h$ does not vanish on $G$ if and only if there exists a {\it thorny vertex} $\vv_i$ and a cut domain
$D_{i_r}\subset D^{\ast}_{\vv_i}$ such that $i_r\in \Lambda^{\ast}_G$, which means that $G\subset D^{\partial}\cup\Omega^c$.
In particular, we have
$$ \tilde{p}_h^{\ast}=0~~~\mbox{on}~~\Omega^{\partial}\backslash (D^{\partial}\cup\Omega^c).$$

Let $G\subset D^{\partial}\cup\Omega^c$. By the definition of $\tilde{p}_{h}^{\ast}$ and (\ref{stability.new}), we can deduce that
$$ \|\tilde{p}_{h}^{\ast}\|_{1, G}\stl\log(1/h)\sum_{i_r\in\Lambda^{\ast}_G}\|p_{h}^{\ast}\|_{1,D_{i_r}}\leq\log(1/h)\sum_{i_r\in\Lambda^{\ast}_G}\|p_{h,\vv_i}^{\ast}\|_{1,\Omega_{i_r}}. $$
This, together with (\ref{stab-5.1new}), gives
\ee
\|\tilde{p}^{\ast}_h\|_{1,G}\stl\log^2(1/h)\sum_{i_r\in\Lambda^{\ast}_G}\|\v_h\|_{\c, \Omega_{i_r}}.\label{stab-5.1newnew}
\e
From the definitions of $D^{\ast}$, $D^{\partial}$ and $\Omega^c$, we have $\al_G\leq\al_{i_{r}}$, and we further obtain $\be_G\leq \mu\be_{i_{r}}$
by {\bf Condition A} in Definition \ref{class}. It follows from (\ref{stab-5.1newnew}) that
$$ \be^{{1\over 2}}_G\|\tilde{p}^{\ast}_h\|_{1,G}\stlm\log^2(1/h)\sum_{i_r\in\Lambda^{\ast}_G}\be^{{1\over 2}}_{i_r}\|\v_h\|_{\c,\Omega_{i_r}}. $$
This, together with {\bf Condition C} in the definition of the class $\mathbb{P}(\gamma)$, leads to
$$ \|\tilde{p}^{\ast}_h\|_{H_{\be}(\Omega^{\partial})}\stlm\log^2(1/h)\|\v_h\|_{H^{\ast}(\c,\Omega^{\ast})}. $$
Note that $\tilde{p}^{\ast}_h=p^{\ast}_h$ on $D^{\ast}$. Thus, using (\ref{stab-5.1new}) again,
\ee
\|\tilde{p}^{\ast}_h\|_{H_{\be}(\Omega)}\stlm\log^2(1/h)\|\v_h\|_{H^{\ast}(\c,\Omega)}. \label{stab-final-5.1}
\e
Similarly, we have
\ee \|\tilde{\w}^{\ast}_{h}\|_{H^1_{\ast}(\Omega)}+h^{-1}\|\tilde{\bf R}^{\ast}_{h}\|_{L^2_{\omega}(\Omega)}
\stlm \log^2(1/h)\|\v_h\|_{H^{\ast}(\c, \Omega)}.\label{stab-final-5.2}
\e

It follows from (\ref{decom-5.6new}) that
\beqx
\|\v^{\partial}_h\|_{H^{\ast}(\c,\Omega^{\partial})}&\stlm&\|\v_h\|_{H^{\ast}(\c,\Omega^{\partial})}+\|\tilde{p}^{\ast}_h\|_{H^1_{\beta}(\Omega)}\cr
&+&\|\r_h\tilde{\w}_{h}^{\ast}\|_{H^{\ast}(\c,\Omega)}+\|\tilde{\bf R}_{h}^{\ast}\|_{H^{\ast}(\c,\Omega)}.
\eqx
Substituting this into (\ref{stab-5.3new}) and using (\ref{stab-final-5.1})-(\ref{stab-final-5.2}), yields
\ee
\|\tilde{p}^{\partial}_h\|_{H^1_{\beta}(\Omega)}\stlm
  \rho^m(h)\log^2(1/h)\|\v_h\|_{H^{\ast}(\c, \Omega)}.\label{stab-final-5.5}
\e
Similarly, we have
\ee
\|\tilde{\w}^{\partial}_{h}\|_{H^1_{\ast}(\Omega)}+h^{-1}\|\tilde{\bf R}^{\partial}_{h}\|_{L^2_{\omega}(\Omega)}
\stlm\rho^m(h)\log^2(1/h)\|\v_h\|_{H^{\ast}(\c, \Omega)}.\label{stab-final-5.6}
\e

Finally, using the definitions of $p_h$, $\w_h$ and ${\bf R}_h$ given in {\bf Step 2}, and combining (\ref{stab-final-5.1})-(\ref{stab-final-5.6}),
we obtain the desired stability estimates (\ref{stab3.newnew1}) and (\ref{stab3.newnew2}).

\hfill $\Box$

\begin{remark}
If there are no {\it thorny vertices}, we have $V_h^{\ast}(\Omega)=V_h(\Omega)$
and $D^{\ast}=\emptyset$. In this case, {\bf Step 1} in the proof is unnecessary and the above logarithmic factor can be dropped.
\end{remark}

\no{\bf Acknowledgments}. The author would like to thank the anonymous reviewer, who gave many insightful suggestions to improve the presentation of this paper.
Moreover, the author wishes to thank Professor Jinchao Xu for suggesting the topic of the
work and giving many knowledgeable comments to improve this article. In particular, Professor Xu alerted
the author to two interesting papers on the {\it quasi-monotonicity assumption}.

\end{document}